\documentclass[12pt,reqno]{amsart}
\usepackage{xifthen}
\usepackage{subfig}
\usepackage[normalem]{ulem}
\usepackage{pkg}
\usepackage{dsfont}
\newcommand{\jansonluczakrucinski}{MR1782847}
\newcommand{\mm}{MR2518205} 

\newcommand{\mezparvir}{MR1026102}
\newcommand{\talagrand}{MR2731561}
\newcommand{\ww}{MR1424469}
\newcommand{\aizenmansimsstarr}{PhysRevB.68.214403}

\newcommand{\bayatigamarniktetali}{MR2743259}
\newcommand{\bbck}{MR1757955}

\newcommand{\bsrecurrence}{MR1873300}

\newcommand{\chvatal}{MR534946}

\newcommand{\dembokaganshepp}{MR1828509}

\newcommand{\dmi}{MR2650042}
\newcommand{\dmsurvey}{MR2643563}
\newcommand{\dommersetal}{MR2733399}

\newcommand{\mms}{springerlink:10.1007/s00440-010-0315-6}

\newcommand{\pnas}{MR2317690}

\newcommand{\mcdiarmid}{MR1036755}

\newcommand{\bcklleftright}{arXiv:1002.0115}
\newcommand{\cdgspotts}{arXiv:1106.4714v3}
\newcommand{\dms}{arXiv:1110.4821}
\newcommand{\slysun}{arXiv:1203.2602}

\newcommand{\rt}{o}
\newcommand{\edge}{\mathrm{\textup{e}}}
\newcommand{\treereg}{\T}
\newcommand{\lpc}{\to_{\mathrm{\textit{loc}}}}
\newcommand{\config}{\mathrm{\textup{cm}}}
\newcommand{\nt}{\ze}

\newcommand{\kd}[2]{\Ind{#1=#2}}
\newcommand{\perm}{\mathrm{\textup{p}}}
\newcommand{\spins}{\mathscr{X}}
\newcommand{\vx}{\mathrm{\textup{vx}}}
\newcommand{\psimin}{\psi_{\min}}
\newcommand{\psimax}{\psi_{\max}}

\newcommand{\es}{\vpi} 
\newcommand{\numin}{\nu_-}
\newcommand{\esmin}{\es_-}
\newcommand{\pimin}{\pi_-}

\newcommand{\bpb}{\mathrm{\textsc{bp}}}
\newcommand{\free}{\mathrm{\textup{f}}}
\newcommand{\plus}{\mathrm{\textup{m}}}
\newcommand{\forp}{{\mathrm{\textup{s}}}} 

\newcommand{\sym}{\mathrm{\textup{sym}}}
\newcommand{\rdc}{\mathbf{R}}
\newcommand{\err}{\mathrm{\textup{err}}}
\newcommand{\zh}{z_h}
\newcommand{\hh}{\mathds{h}}

\newcommand{\splx}{\triangle}
\newcommand{\splxstar}{\splx^\star}
\newcommand{\splxpr}{\splx_\edge}
\newcommand{\splxprpm}{\splx^\pm_\edge}
\newcommand{\splxbal}{\bar\splx}
\newcommand{\splxfk}[1]{\splxbal[#1]}

\newcommand{\usi}{\underline{\smash{\si}}}
\newcommand{\upsi}{\underline{\smash{\psi}}}
\newcommand{\uh}{\underline{\smash{h}}}
\newcommand{\ueta}{\underline{\smash{\eta}}}
\newcommand{\ub}{\underline{\smash{b}}}

\newcommand{\vpsi}{{\bar\psi}}
\newcommand{\vxi}{{\bar\xi}}
\newcommand{\vh}{{\bar h}}

\newcommand{\vde}{\bar\de}

\newcommand{\bPhi}{{\bm{\Phi}}}
\newcommand{\bh}{\bm{h}}
\newcommand{\bbh}[2]{\bm{h}_{#1 #2}}
\newcommand{\bde}{{\bm{\delta}}}
\newcommand{\bbde}[2]{{\bm{\delta}_{#1 #2}}}
\newcommand{\bph}{{\bm{\varphi}}}
\newcommand{\bbph}[2]{{\bm{\varphi}_{#1 #2}}}
\newcommand{\bA}{\mathbf{A}}
\newcommand{\bzh}{\bm{z}_h}

\newcommand{\CC}{\mathfrak{C}}
\newcommand{\MM}{\mathfrak{M}}

\newcommand{\ratio}{\mathscr{R}}

\newcommand{\BP}{\mathrm{\textup{\textsf{BP}}}}
\newcommand{\pp}{\mathrm{\textup{\textsf{p}}}}
\newcommand{\QQ}{\mathrm{\textup{\textsf{Q}}}}
\newcommand{\oo}{\mathrm{\textup{\textsf{o}}}}
\newcommand{\bpff}{\mathrm{\textup{\textsf{F}}}}
\newcommand{\rr}{\mathrm{\textup{\textsf{r}}}}

\title[Replica symmetry for Potts models on $d$-regular graphs]{The replica symmetric solution \\ for Potts models on $d$-regular graphs}
\author[A.\ Dembo]{Amir Dembo$^*$}
\author[A.\ Montanari]{Andrea Montanari$^\dagger$}
\author[A.\ Sly]{Allan Sly$^\ddagger$}
\author[N.\ Sun]{Nike Sun$^\S$}
\address{$^*$Department of Mathematics, Stanford University \newline\indent Building 380, Sloan Hall, Stanford, California 94305}
\address{$^\dagger$Department of Electrical Engineering, Stanford University \newline\indent 350 Serra Mall, Stanford, California 94305}
\address{$^\ddagger$Department of Statistics, University of California, Berkeley \newline\indent Evans Hall, Berkeley, California 94720}
\address{$^{*\dagger\S}$Department of Statistics, Stanford University \newline\indent Sequoia Hall, 390 Serra Mall, Stanford, California 94305}
\date{\today}
\thanks{Research partially supported by NSF grants $^{*\dagger\S}$DMS-1106627 and $^\dagger$CCF-0743978, $^\ddagger$Alfred P.\ Sloan Research Fellowship, and $^\S$Department of Defense NDSEG Fellowship.}

\begin{document}

\maketitle

\begin{abstract}
We provide an explicit formula for the limiting
free energy density (log-partition function divided by the number of vertices) for ferromagnetic Potts models on uniformly sparse graph sequences converging locally to the $d$-regular tree for $d$ even, covering \emph{all} temperature regimes. This formula coincides with the Bethe free energy functional evaluated at a suitable fixed point of the belief propagation recursion on the $d$-regular tree, the so-called replica symmetric solution. For uniformly random $d$-regular graphs we further show that the replica symmetric Bethe formula is an upper bound for the asymptotic free energy for \emph{any} model with permissive interactions.
\end{abstract}

\section{Introduction}\label{s:intro}

Let $G=(V,E)$ be a finite undirected graph, and $\spins$ a finite
alphabet of \emph{spins}. A \emph{factor model} on $G$ is a
probability measure on the space of \emph{(spin) configurations}
$\usi\in\spins^V$ of  the form
\beq\label{e:fm}
\nu^{\upsi}_G(\usi)
\equiv \f{1}{Z_G(\upsi)}
	\prod_{(ij)\in E} \psi(\si_i,\si_j)
	\prod_{i\in V} \vpsi(\si_i),
\eeq
where $\psi$ is a symmetric function $\spins^2\to\R_{\ge0}$, $\vpsi$ is a positive function $\spins\to\R_{>0}$, and $Z_G(\upsi)\equiv Z_G$ is the normalizing constant, called the \emph{partition function} (with its logarithm called the \emph{free~energy}). The pair $\upsi\equiv(\psi,\vpsi)$ is called a \emph{specification} for the factor model~\eqref{e:fm}. We assume the specification is \emph{permissive}, that is, there exists $\si^\perm\in\spins$ with $\min_\si\psi(\si,\si^\perm)>0$.

A primary example we consider in this paper is the \emph{$q$-state Potts model} on $G$ with inverse temperature $\be$ and magnetic field $B$, given by specification
\beq\label{e:potts}
\psi(\si,\si') = e^{\be\kd{\si}{\si'}},\quad
\vpsi(\si) = e^{B\kd{\si}{1}},\quad
\spins=[q]\equiv\set{1,\ldots,q}.
\eeq
We write $\nu^{\be,B}_G$ for the corresponding measure on $[q]^V$. The model is said to be \emph{ferromagnetic} if $\be\ge0$, and \emph{anti-ferromagnetic} otherwise.

In this paper we study the asymptotics of the free energy for factor models \eqref{e:fm} on graph sequences $G_n=(V_n,E_n)$ converging locally to the $d$-regular tree $\treereg_d$ ($d\ge3$) in the sense of Benjamini--Schramm~\cite{\bsrecurrence} (see Defn.~\ref{d:lwc}). This class includes in particular any sequence of $d$-regular graphs with girth (minimal cycle length) diverging to infinity.

The study of statistical mechanics on regular trees has a long history, initiated by Bethe~\cite{Bethe}. While tree graphs do not capture the finite-dimensional structure of actual physical systems, models on trees are often amenable to exact analysis. Further, it is often argued that they are a good approximation to models on the lattice $\Z^d$ for large $d$ or for long interaction range~\cite{Weiss,Abou-Chacra,Chalupa,ThoulessSG}. According to this expectation, models on trees provide a flexible and well-defined approach for investigating \emph{mean-field theory} (i.e.~the behavior of statistical mechanics models in high dimensions).

While this expectation proves to be correct in a number of examples, it has recently become clear that, in many cases, models on trees fail to capture the ``correct'' mean-field behavior. Spin glasses provide an important example of this phenomenon: a fairly
natural class of spin glasses on trees was introduced by
Thouless~\cite{ThoulessSG} and further characterized by Chayes et al.~\cite{ChayesSG}. However, the thermodynamic behavior observed there is very different from the widely accepted mean-field theory of spin glasses, as obtained from analysis of the Sherrington--Kirkpatrick ($\acr{sk}$) model~\cite{\mezparvir,\talagrand}. In particular, the low-temperature phase of the tree models defined in~\cite{ThoulessSG} does not exhibit replica symmetry breaking (in contrast with $\acr{sk}$). A similar discrepancy was observed in the case of Anderson localization by Aizenman--Warzel~\cite{AizenmanWarzel}.

In the case of spin  glasses, M\'ezard--Parisi~\cite{MezardParisiBethe} argued that this difference arises because of a particular feature of tree graphs: in the subgraph induced by the first $\ell$ levels of the regular tree, the leaves constitute a non-vanishing fraction of the vertices as $\ell\to\infty$. They suggested that mean-field theory ought instead to be defined by considering graphs that are not themselves trees, but ``look like regular trees'' in the neighborhood of a typical vertex (which fails for the depth-$\ell$ subtree of the regular tree) --- the canonical example being the (uniformly) random $d$-regular graph ensemble. This approach allows to reconcile discrepancies in several known cases. In particular, spin glasses on random regular graphs are expected to exhibit replica symmetry breaking with features analogous to the $\acr{sk}$ model (see~\cite{MezardParisiBethe} and \cite[Ch.~17]{\mm}).

Let us also mention that the study of statistical mechanics models on locally tree-like graphs has attracted renewed interest because of the connection with random combinatorial problems, such as $k$-$\acr{sat}$ and graph coloring. Statistical physicists were indeed able to compute threshold locations for these models by analyzing suitable Gibbs measures on locally tree-like structures~\cite{MezParZec_science,\pnas,\mm}. Rigorous verification of these predictions is an outstanding mathematical challenge.

In this paper we consider the existence and value of the \emph{free energy density} (asymptotic free energy per spin)
\beq\label{e:lim}
\phi
\equiv \lim_{n\to\infty}\phi_n
\equiv \lim_{n\to\infty} \f{1}{n}\E_n[\log Z_n],
\quad Z_n\equiv Z_{G_n}(\upsi),
\eeq
for $G_n$ a  (possibly random) graph sequence converging locally to the regular tree and $\E_n$ expectation over the law of $G_n$. For Ising (specification~\eqref{e:potts} with $q=2$) models in the ferromagnetic regime, for \emph{any} graph sequence with uniformly integrable average degree converging locally to a (possibly random) tree, the free energy density \eqref{e:lim} exists and depends only on the limiting tree~\cite{\dmi,\dommersetal,\dms}. The computation of $\phi$ allows to compute various limits of interest with respect to the $\nu_{G_n}$, as done for example in \cite{\mms,\dommersetal}. Proving existence of the free energy density for $q\ge2$ and general specification $\upsi$ poses several challenges:\footnote{Existence of \eqref{e:lim} for general $\upsi$ is equivalent to right convergence of $G_n$ in the language of~\cite{\bcklleftright}.}

\smallskip
\noindent
{\it 1.} There are examples in which the free energy density \eqref{e:lim} depends not only on the limiting tree but also on the particular graph sequence. For example, in the anti-ferromagnetic Ising model at sufficiently low temperature (sufficiently negative $\be$), it is not difficult to show that the free energy per spin on random $d$-regular graphs is asymptotically lower than on random \emph{bipartite} $d$-regular graphs. As a consequence local weak convergence is not in full generality a sufficient condition for existence of the limit \eqref{e:lim}.

\smallskip
\noindent
{\it 2.} Statistical physicists have put forth a number of conjectures (corresponding to different models or regimes) on the free energy density \eqref{e:lim} (see e.g.~\cite{\mezparvir,\mm}). This analysis generally imposes a probability distribution on the graph $G_n$ which is suitable for calculations, typically the Erd\doubleacute{o}s-Renyi or configuration models. Ensuing rigorous work has also focused on the same random graph ensembles (see e.g.\ \cite{\talagrand}) rather than understanding which graph sequences in general have a limit \eqref{e:lim}. In this paper we focus instead on individual graph sequences.

\smallskip

Characterizing the limit for ensembles of uniformly random graphs is already beyond current techniques for many factor models \eqref{e:fm}. Achieving the same goal for general locally tree-like graph sequences is all the more difficult, and requires to go beyond what is known from physics methods. A simple example is provided again by the anti-ferromagnetic Ising model: existence of the limit can be proved by a combinatorial interpolation~\cite{\bayatigamarniktetali}, but even a heuristic prediction of the value is unavailable.

In contrast, as mentioned above the free energy density for the ferromagnetic Ising model on locally tree-like graphs exists and can be computed. Its value is given by the \emph{Bethe prediction} $\Phi$, which is expressed in terms of a certain fixed point of a distributional recursion (given in \eqref{e:bp} in the $d$-regular setting; see~\cite{\dmsurvey,\dms} for the general case). This result was proved in the case of Galton--Watson limiting trees in~\cite{\dmi} via an interpolation scheme. In~\cite{\dms} a generalized scheme was developed which gave the result for Ising on general limiting trees. The method was applied also to show $\liminf_n\phi_n\ge\Phi$ in the ferromagnetic $q$-Potts model with $q>2$, but could only pin down $\phi=\Phi$ in limited regimes of $(\be,B)$. The difficulty of the Potts model with $q>2$ may be understood as follows: by a monotonicity argument, the local weak limit of Potts measures on $G_n$ is sandwiched between the free and maximally $1$-biased Gibbs measures on the limiting tree. When $q=2$ these measures coincide for any $\be\ge0,B>0$,\footnote{Equivalently, there is only one Gibbs measure $\nu$ on $\treereg_d$ that satisfies the following properties: (i) $\nu$ is invariant under automorphisms of $\treereg_d$; (ii) $\nu$ is a Markov chain on $\treereg_d$; and (iii) $\nu(\sigma_\rt=+1)>0$.} but when $q>2$ the measures disagree in certain regimes of $(\be,B)$. This corresponds to the appearance of ``multiple stable fixed points'' in the distributional recursion \eqref{e:bp} as soon as $q>2$.

In this paper we establish the existence of the free energy density \eqref{e:lim}, and provide an explicit expression for its value, on graphs converging to regular trees of even degree and for all $q>2$. Let us mention that the statistical physics folklore prescribes that the distributional fixed point with the \emph{highest Bethe free energy density} should be selected. However, in the  physics literature this is justified only via analogy with other models, without providing arguments which apply to locally tree-like graphs. Our result is the first rigorous verification of this variational principle in a non-trivial example for locally tree-like graphs.

A different variational principle was proved in~\cite{GuerraRSB,\aizenmansimsstarr} for mean-field spin glass models, but in that case the free energy density needs to be \emph{minimized}. This difference is typically attributed by physicists to the difference between ferromagnetic and spin glass models; it remains an outstanding challenge to understand these two variational principles within a common framework. In the context of models on sparse graphs, Contucci et al.\ \cite{\cdgspotts} recently proved that the variational principle of~\cite{GuerraRSB,\aizenmansimsstarr} provides a bound on the free energy of \emph{anti-ferromagnetic} Potts models, which was proved to be tight at high temperature.

The rest of the paper is organized as follows: in the remainder of this introductory section we review the definitions of local convergence and the Bethe prediction and formally state our results, which we divide into two categories: in \S\ref{s:unif} we study the Bethe prediction on the \emph{uniformly} random $d$-regular graph ensemble. In \S\ref{s:decomp}-\ref{s:potts} we prove results in the more general setting of graphs converging locally to the $d$-regular tree. In each case we first consider general specifications $\upsi$ before specializing the the Potts specification \eqref{e:potts}.

\subsection{Local convergence}

If $G$ is any graph and $U$ any subgraph, we write $\pd U$ for the external boundary of $U$ in $G$ (the set of vertices in $G$ adjacent to but not contained in $U$). For any vertex $v$ of $G$, we let $D_v\equiv\abs{\pd v}$ denote its degree, and write $\ball{t}{v}$ for the subgraph induced by the vertices of $G$ at graph distance at most $t$ from $v$. Fix $d$ throughout and let $\treereg_d\equiv(\treereg_d,\rt)$ denote the $d$-regular tree rooted at $\rt$, with $\treereg^t_d\equiv\ball{t}{\rt}$ the subtree of depth $t$.

\bdfn\label{d:lwc}
For $G=(V,E)$ finite undirected, let $\nt_t(G)\equiv\abs{V}^{-1}\abs{\set{v\in V : \ball{t}{v}\not\cong\treereg^t_d}}$ where $\cong$ denotes graph isomorphism. The sequence of (random) graphs $G_n=(V_n=[n],E_n)$ is said to \emph{converge locally to the $d$-regular tree $\treereg_d$} if for all $t\ge0$, $\nt_t(G_n)\to0$ in probability as $n\to\infty$.
\edfn

For $G=(V,E)$ let $I_G$ denote a vertex chosen uniformly at random from $V$, and write $I_n\equiv I_{G_n}$. From now on let $\P_n$ denote the joint law of $(G_n,I_n)$, and $\E_n$ the expectation with respect to $\P_n$. An equivalent definition of the local convergence of $G_n$ to $\treereg_d$ is that $\lim_{n\to\infty}\P_n(\ball{t}{I_n}\cong\treereg^t_d)=1$ for all $t\ge0$.

\bdfn
The sequence $G_n$ is \emph{uniformly sparse} if the random variables $D_{I_n}$ are uniformly integrable, that is, if $$\lim_{L\to\infty}\limsup_{n\to\infty}\E_n[D_{I_n} \Ind{ D_{I_n}\ge L }]=0.$$
We assume throughout that $G_n$ ($n\ge1$) is a uniformly sparse graph sequence converging locally to the $d$-regular tree $\treereg_d$. This setting is hereafter denoted $G_n\lpc\treereg_d$.
\edfn

\subsection{The Bethe prediction}

\subsubsection{Definition in $d$-regular setting}

We now describe the Bethe free energy prediction in the special setting of $d$-regular trees; for a more general description see~\cite{\dmsurvey,\dms}. We write $\si$ for elements of the finite alphabet $\spins$ of spins, and $\usi$ for vectors with entries in $\spins$; supposing first that $d$ is \emph{even}, for $\usi\in\spins^d$ let
$$\Om^\vx(\usi)\equiv\sum_\si\vpsi(\si)\prod_{j=1}^d\psi(\si,\si_j),\quad
\Om^\edge(\usi)\equiv\prod_{j=1}^{d/2} \psi(\si_{2j-1},\si_{2j})$$
Let $\splx_{\spins^k}$ denote the $(\abs{\spins}^k-1)$-dimensional simplex of probability measures on $\spins^k$.
If $\hh$ is a finite measure on $\spins^k$ (any $k$) and $g$ is any function on $\spins^k$, then $\angl{g}_\hh$ denotes the integral of $g$ with respect to $\hh$. For $\hh\in\splx_{\spins^d}$ let
\beq\label{e:Psi}
\Psi^\vx(\hh)\equiv
	\angl{\Om^\vx}_\hh,\quad
\Psi^\edge(\hh)\equiv
	\angl{\Om^\edge}_\hh,\quad
\Psi(\hh)\equiv\f{\Psi^\vx(\hh)}{\Psi^\edge(\hh)}.
\eeq
With a slight abuse of notation we write $\hh\equiv\uh\equiv(h^1,\ldots,h^d)$ with $h^j\in\splx\equiv\splx_\spins$ to indicate that $\hh$ is the product measure $\hh(\usi)=\prod_{j=1}^d h^j_{\si_j}$. The \emph{Bethe free energy functional} is then defined for $h\in\splx$ by
\begin{align*}
\Phi(h)
&\equiv\Phi^\vx(h)-\Phi^\edge(h)
\equiv\log\Psi^\vx(h,\ldots,h)-\log\Psi^\edge(h,\ldots,h)\\
&=
\log\bigg\{ \sum_\si\vpsi(\si)
	\Big( \sum_{\si'}\psi(\si,\si')h_{\si'} \Big)^d \bigg\}
-\f{d}{2}\log\bigg\{ \sum_{\si,\si'}
	\psi(\si,\si')h_\si h_{\si'} \bigg\};
\end{align*}
this definition clearly extends to $d$ odd. The \emph{Bethe prediction} is that the asymptotic free energy $\phi$ of \eqref{e:lim} exists and equals
\beq\label{e:Phi}
\Phi\equiv\sup_{h\in\splxstar}\Phi(h),
\eeq
where $\splxstar$ denotes the set of fixed points in $\splx$ of the \emph{belief~propagation} or \emph{Bethe~recursion} $\BP:\splx\to\splx$, defined by
\beq\label{e:bp}
(\BP h)(\si)
\equiv \f{1}{\zh}
\vpsi(\si)
	\Big( \sum_{\si'}\psi(\si,\si')h_{\si'} \Big)^{d-1},
\eeq
with $\zh$ the normalizing constant. For permissive $\upsi$, any fixed point $h\in\splxstar$ must belong to the interior of $\splx$ (i.e.\ $\min_\si h_\si>0$). An interior point $h$ of $\splx$ belongs to $\splxstar$ if and only if
\beq\label{e:bpval}
\f{\vpsi(\si)}{h_{\si}}
	\Big(\sum_{\si'}\psi(\si,\si') h_{\si'}\Big)^{d-1}
	=\zh
	\quad\forall\si\in\spins.
\eeq
In this case, writing $\bzh\equiv\angl{\psi}_{h\otimes h}\equiv\sum_{\si,\si'} \psi(\si,\si')h_\si h_{\si'}$,
\beq\label{e:Phi.simp}
\Phi^\vx(h)=\log\zh+\log\bzh,\quad
\Phi^\edge(h)=\f{d}{2}\log\bzh.
\eeq
Fixed points $h\in\splxstar$ correspond to ``Bethe Gibbs measures,'' suitable candidates for the local weak limit of $\nu_{G_n}$ (see e.g.\ \cite[Rmk.~1.12]{\dms}).

\subsubsection{Bethe variational principle}

Let $\splxpr$ denote the (compact) set of symmetric probability measures $\bh$ on $\spins^2$, with one-point marginals denoted by $\vh$. Then $\splxstar$ embeds into $\splxpr$ via the relation
\beq\label{e:bij}
\bbh{\si}{\si'}
= (h\otimes_\psi h)_{\si\si'}
\equiv \f{\psi(\si,\si')h_\si h_{\si'}}{\bzh}.
\eeq
Let $\xi\equiv\log\psi$ and $\vxi\equiv\log\vpsi$. We then define
\begin{align}
\nonumber \bPhi(\bh)
&\equiv \angl{\vxi}_\vh-(d-1)H(\vh)+\f{d}{2}[\angl{\xi}_{\bh}+H(\bh)]\\
\label{e:Phi.var} &=-\relent{\vh}{\vpsi}
-\f{d}{2} \relent{\bh}{\vh\otimes_\psi \vh}.
\end{align}
In the above and hereafter, for $p,q$ finite non-negative measures on $\spins^k$, $H(p)$ denotes the Shannon entropy $-\sum_x p_x \log p_x$, and $\relent{q}{p}$ denotes the relative entropy $\sum_x q_x \log (q_x/p_x)$ between $q$ and $p$. We take the usual conventions $\log0=-\infty$, $0\log 0=0$ and $0\log(0/0)=0$.

The Bethe prediction has the following variational characterization:

\bppn\label{p:var}
Let $\upsi\equiv(\psi,\vpsi)$ be a permissive specification.
\bnm[(a)]
\item \label{p:var.a} Any interior stationary point $\bh$ of $\bPhi$ corresponds to $h\in\splxstar$ by the bijective relation \eqref{e:bij}. Any local maximizer $\bh$ of $\bPhi$ is an interior point of $\splxpr$, so
\beq\label{e:var}
\Phi=\sup_{\bh\in\splxpr}\bPhi(\bh).
\eeq
\item \label{p:var.b} An interior stationary point $\bh$ of $\bPhi$ is a local maximizer if and only if, for $(X,Y)$ having (exchangeable) law $\bh$,
\beq\label{e:corr}
\rh_{XY}\equiv \sup\bigg\{
\f{\Var \E[\bbph{X}{Y}  \giv X]}{\Var \bbph{X}{Y} } :
	\bph\not\equiv0, \bbph{\si}{\si'}=\bbph{\si'}{\si}
\bigg\}
\le \f{d}{2(d-1)}.
\eeq
\enm

\bpf
\eqref{p:var.a} Follows from \cite[Thm.~1.16]{\dms} (using compactness of $\splxpr$).

\medskip\noindent\eqref{p:var.b}
Let $\splxprpm$ denote the set of functions $\bde:\spins^2\to\R$ satisfying
$$\bbde{\si}{\si'}=\bbde{\si'}{\si},\quad
\sum_{\si,\si'}\bbde{\si}{\si'}=0,\quad\text{and}\quad
\sum_{\si,\si'}\bbde{\si}{\si'}^2=1.$$
It was shown in \cite[Propn.~3.4]{\dms} that an interior stationary point $\bh$ of $\bPhi$ is a local maximizer if and only if
\beq\label{e:locmax.delta}
\left.4\,\pd_\eta^2 \Phi_\mu(\bh+\eta\bde)\right|_{\eta=0}
=2(d-1) \angl{(\vde/\vh)^2}_{\vh}-d\angl{(\bde/\bh)^2}_{\bh}\le0
\quad\forall\bde\in\splxprpm.
\eeq
The condition \eqref{e:corr} follows by taking $\bph=(\bh+\bde)/\bh$ and rearranging.
\epf
\eppn

\brmk
The ``symmetric correlation coefficient'' $\rh_{XY}$ measures dependence of the exchangeable pair $(X,Y)$.\footnote{Note $\rh_{XY}$ is not the classical correlation coefficient between $\si(X), \si(Y)$ (see~e.g.~\cite{\dembokaganshepp} and references therein).} By the classical variance decomposition, $\rh_{XY}\in[0,1]$ with $\rh_{XY}=1$ if and only if $Y=f(X)$ for some deterministic function $f$ (which by exchangeability must be involutive). If $X$ and $Y$ are independent, it is easily seen from Hoeffding's decomposition $\bbph{X}{Y} \equiv \wt\bph_{XY} +\E[\bbph{X}{Y} \giv X]+\E[\bbph{X}{Y} \giv Y]-\E[\bbph{X}{Y} ]$ (with $\E[\wt\bph_{XY} \giv X]=0$) that $\rh_{XY}=1/2$ with supremum achieved by $\bbph{\si}{\si'}$ of form $\vph_\si +\vph_{\si'}$. We do not know of an argument to show $\rh_{XY}\ge1/2$ for any exchangeable $(X,Y)$.
\ermk

\subsection{Results for uniformly random $d$-regular graphs}
\label{ss:intro.unif}

\subsubsection{Expectation of the partition function}
\label{sss:intro.unif.exp}

For $dn$ even let $\cM_{d,n}$ be the space of perfect matchings of $[dn]$, and for $\gm\in\cM_{d,n}$ let $G[\gm]$ be the (multi-)graph (i.e.\ with multi-edges and self-loops permitted) on vertex set $[n]$ defined by $\gm$ through the projection $[dn]\to[n]$ taking $i'\in[dn]$ to its representative modulo $n$ in $[n]$. The \emph{configuration model} is the probability measure $\P^\config_{d,n}$ on $\cG^\config_{d,n}\equiv\set{G[\gm]:\gm\in\cM_{d,n}}$ induced by the uniform measure on matchings $\cM_{d,n}$. The measure $\P^\config_{d,n}$ conditioned on the set $\cG_{d,n}$ of \emph{simple} graphs is simply the uniform measure $\P_{d,n}$ on the $d$-regular graphs on $[n]$. We write $\phi^\config_{d,n}\equiv n^{-1}\E^\config_{d,n}[\log Z_n]$ where $\E^\config_{d,n}$ denotes expectation under $\P^\config_{d,n}$.

\bThm\label{t:fm.unif}
For any permissive specification $\upsi\equiv(\psi,\vpsi)$,
\beq\label{e:fm.unif}
\phi^\config
\equiv\limsup_{n\to\infty}
	\phi^\config_{d,n}
\le \lim_{n\to\infty}
	\f1n
	\log \E^\config_{d,n}[Z_n]
=\Phi
\eeq
\eThm

If $G_-$ is any (multi-)graph on $[n]$ with maximum degree at most $d$ and $G$ is formed by adjoining a new vertex $i$ to $d$ or fewer vertices in $G_-$, then
\begin{align}
\nonumber
\f{Z_G}{Z_{G_-}}
&=\sum_{\si_i} \vpsi(\si_i)
	\sum_{\usi_{\pd i}}
	\prod_{j\in\pd i}\psi(\si_i,\si_j)
	\nu_{G_-}(\usi_{\pd i})
\le\psimax^{d+1},\quad
\psimax
	\equiv\max_{\si,\si'}\,[\vpsi(\si)\vee\psi(\si,\si')],\\
\label{e:perm.bd}
\f{Z_G}{Z_{G_-}}
&\ge \vpsi(\si^\perm)
	\sum_{\usi_{\pd i}}
	\prod_{j\in\pd i}\psi(\si^\perm,\si_j)
	\nu_{G_-}(\usi_{\pd i})
\ge\psimin^{d+1},\quad
\psimin\equiv\min_\si\,[\vpsi(\si) \wedge \psi(\si,\si^\perm)],
\end{align}
so $\abs{\log Z_G-\log Z_{G_-}}$ is uniformly bounded by a constant depending only on $d,q,\upsi$. Consequently, if $\P_n$ is any probability measure on (multi-)graphs $G_n=(V_n\equiv[n],E_n)$ with maximum degree at most $d$, then the Azuma--Hoeffding bound applied to the vertex-revealing martingale gives $\P_n[\abs{n^{-1}\log Z_n-\phi_n}\ge \ep] \le e^{-c\ep^2 n}$ for some constant $c\equiv c(d,q,\upsi)>0$. On the other hand, even if asymptotically $\phi^\config_{d,n}<\Phi$, there is a substantial subclass of graphs in $\cG^\config_{d,n}$ with free energy close to $\Phi$: since $n^{-1}\log Z_n$ is uniformly bounded over such graphs by a constant $C\equiv C(d,q,\upsi)<\infty$, it follows from \eqref{e:fm.unif} that $p_n=\P^\config_{d,n}[n^{-1}\log Z_n\ge\Phi-\ep]$ satisfies
$$\Phi\le
\liminf_{n\to\infty}\f1n
	\log\big(
	p_n e^{C n}
	+[1-p_n] e^{(\Phi-\ep)n}
	\big),$$
thus $p_n\ge e^{-[C-\Phi+\ep/2]n}$ for large $n$. Since \beq\label{e:config.err}
\abs{\cG_{d,n}}\asymp\abs{\cG^\config_{d,n}}\asymp(nd-1)!!\sim\sqrt{2} (nd/e)^{nd/2}
\eeq
(see e.g.~\cite[Ch.~9]{\jansonluczakrucinski}), this means that the set of graphs in $\cG_{d,n}$ with free energy at least $\Phi-\ep$ also grows like $e^{(d/2)n\log n(1-o(1))}$. In fact, we have the following

\bcor\label{c:fm.unif}
If $\phi^\config<\Phi$, then there exist $x,C_x>0$ with
\beq\label{e:unif.tree.above}
\P^\config_{d,n}[n^{-1}\log Z_n \ge\Phi+x] \ge e^{-C_x n},
\eeq
and a sequence of graphs $G_n\in\cG^\config_{d,n}$, $G_n\lpc\treereg_d$ with $\liminf_{n\to\infty} n^{-1}\log Z_n \ge\Phi+x$.
\ecor

Although the proof of the corollary is straightforward, we highlight it here because it demonstrates that if the uniformly random ensemble has free energy $\phi_n$ strictly below the replica symmetric solution $\Phi$ --- as is expected to happen in replica symmetry breaking regimes --- then there is a breaking of homogeneity in the graph space $\cG^\config_{d,n}$ as well, with a large subclass of graphs achieving free energy strictly above $\Phi$. An interesting open question is whether the maximal asymptotic free energy is achieved by random bipartite graphs, as is known to be the case in two-spin models \cite{\slysun}.

\subsubsection{The Potts Bethe prediction}

Surprisingly, another consequence of Thm.~\ref{t:fm.unif} is the following solution to the optimization problems \eqref{e:Phi} and \eqref{e:Phi.var} for the ferromagnetic Potts model. Let $h^\free$ denote the limit of successive iterations of $\BP$ starting from the uniform probability measure on $[q]$, and let $h^\plus$ denote the limit of successive iterations of $\BP$ starting from the probability measure on $[q]$ supported on spin $1$.

\bThm\label{t:potts.sup}
For the Potts model \eqref{e:potts} with $\be,B\ge0$, $\Phi=\Phi(h^\free)\vee\Phi(h^\plus)$, and if $B>0$ then this is strictly greater than $\Phi(h)$ for any $h\in\splxstar\setminus\set{h^\free,h^\plus}$.
\eThm

We supplement Thm.~\ref{t:potts.sup} by a classification of stationary points of $\bPhi$ (equivalently, via \eqref{e:bij}, solutions of the Potts Bethe recursion) as well as a study of which stationary points can be local maximizers. The motivation for considering \emph{local} maximizers of $\bPhi$ --- which after all are irrelevant to the Bethe prediction \eqref{e:var} if they are not global maximizers --- is that we expect these are precisely the fixed points which can be seen in local weak limits of \emph{conditioned} factor models on graph sequences $G_n\lpc\treereg_d$, in the spirit of~\cite{\mms}. That is, when $\bh$ is a local maximizer of $\bPhi$, the factor model restricted to configurations of edge empirical measure close to $\bh$ should converge locally weakly to the (Bethe) Gibbs measure corresponding to $\bh$.

Detailed statements are given in Propns.~\ref{p:potts.sols} and~\ref{p:potts.locmax}. We show in particular that any local maximizer $\bh$ of $\bPhi$ must correspond (via \eqref{e:bij}) to $h$ with $\abs{\set{h_\si :\si\in\spins}}\le3$. On the other hand we show that for $B>0$ small and $\be>0$ large there exist solutions $h=(\QQ,\pp_+,\pp_-,\ldots,\pp_-)$ corresponding (via \eqref{e:bij}) to local maximizers of $\bPhi$.

\subsection{Results for general $d$-regular graphs}
\label{ss:intro.gen}

For $d$ \emph{even}, we establish the Potts Bethe prediction for general graph sequences $G_n\lpc\treereg_d$:

\bThm\label{t:potts}
For the Potts model on $G_n\lpc\treereg_d$ with $d$ even, $\phi(\be,B)=\Phi(\be,B)$ for all $\be,B\ge0$.
\eThm

This theorem will be deduced from the following result for abstract factor models which illustrates a more general principle. We restrict hereafter to \emph{$d$-regular} graph sequences $G_n\lpc\treereg_d$, since in \S\ref{s:decomp} we will show that, for the purposes of computing the free energy, general sequences $G_n\lpc\treereg_d$ can be reduced to the $d$-regular case using the uniform sparsity hypothesis.

For fixed $G=(V,E)$ we let $\E_G$ denote expectation over the uniformly random vertex $I_G\in V$. We define symmetrized versions of $\Psi^\edge,\Psi$ by
$$\Psi^{\edge,\sym}(\hh)
\equiv\f{1}{d!}
	\sum_{\pi\in S_d} \Psi^\edge(\hh^\pi),\quad
\Psi^\sym
\equiv\f{\Psi^\vx}{\Psi^{\edge,\sym}}$$
where $S_d$ denotes the symmetric group on $d$ letters and $\hh^\pi(\usi)
\equiv
\hh(\si_{\pi(1)},\ldots,\si_{\pi(d)})$. Given a measure $\rh$ on $\splx^d$ define the mixture of product measures
\beq\label{e:bar.rho}
\bar\rh\in\splx_{\spins^d},\quad
\bar\rh(\usi)
\equiv\int h^1_{\si_1}\cdots h^d_{\si_d} \,d\rh(\uh),
\eeq
Throughout we write $\oo(t,x)$ for a uniformly bounded function such that
$$\lim_{x\decto0}\oo(t,x)=\oo(t),\quad
\lim_{t\to\infty}\oo(t)=0.$$
The function may change from line to line with the understanding that it can be chosen to depend only on $d,q,\upsi$.

\bThm\label{t:fm}
Suppose $\upsi$ is a permissive specification and $G_n\lpc\treereg_d$, $d$ even. Suppose that for all finite $d$-regular graphs $G$, all $v\in G$, and all $t\in\Z_{\ge0}$, there is a measure $\rh\equiv\rh(G,v,t)$ on $\splx^d$ such that
\bnm[(i)]
\item\label{cond:prod}
For $I=I_G$ and $\numin\equiv\nu_{G\setminus I}(\usi_{\pd I}=\cdot)$, the measure $\bar\rh$ satisfies $\E_G[ \tv{\numin-\bar\rh}]=\oo(t,\nt_t(G))$.

\item\label{cond:ubd}
For all $v,t$ we have $\log\Psi^\sym(\uh)\le\Phi+\oo(t)$ for all $\uh\in\supp\rh(G,v,t)$.
\enm
Then $\limsup_n\phi_n\le\Phi$ for the factor model on $G_n$ specified by $\upsi$.
\eThm

In \S\ref{s:potts} we show that the conditions of the preceding theorem are satisfied in the Potts model with $\be,B\ge0$ and hence $\limsup_{n\to\infty}\phi_n(\be,B)\le\Phi(\be,B)$; Thm.~\ref{t:potts} is then proved as the matching lower bound $\liminf_{n\to\infty}\phi_n(\be,B)\ge\Phi(\be,B)$ was shown in \cite[Thm.~1.10]{\dms}. To give some motivation for the conditions of Thm.~\ref{t:fm}, consider the following randomized operation $\rdc$ on finite $d$-regular graphs $G=(V,E)$:
\begin{quote}
\emph{Operation $\rdc$:} \\
Let $(v_1,\ldots,v_d)$ be an enumeration of the neighbors of a uniformly random vertex $I\equiv I_G$ in $G$. For each $\pi\in S_d$ let $G^\pi$ be the graph formed by adding to $G_-\equiv G\setminus I$ the edges $(v_{\pi(2j-1)},v_{\pi(2j)})$ for $1\le j\le d/2$. Then set $\rdc G\equiv G^\pi$ for $\pi\in\argmin_{\pi'}[\log Z_G-\log Z_{G^{\pi'}}]$.
\end{quote}
Thm.~\ref{t:fm} is then proved (in \S\ref{s:decomp}) by expressing $\E_n[\log Z_n]$ as a telescoping sum over $\E_n[\log Z_{\rdc^j G}-\log Z_{\rdc^{j+1} G}]$, and showing that conditions \eqref{cond:prod} and \eqref{cond:ubd} above imply that each term in the sum is bounded above by $\Phi$: indeed, with $\numin\equiv\nu_{G_-}(\usi_{\pd I}=\cdot)\in\splx_{\spins^d}$ as before, notice that
\beq\label{e:ratio.gen}
\f{Z_G}{Z_{G_-}}=\angl{\Om^\vx(\usi)}_{\numin}
=\Psi^\vx(\numin),\quad
\f{Z_{G^\pi}}{Z_{G_-}}
=\angl{\Om^\edge(\usi)}_{\numin^\pi}
=\Psi^\edge(\numin^\pi).
\eeq
Condition \eqref{cond:prod} gives that these are well approximated by
\beq\label{e:ratio.prod}
\Psi^\vx(\bar\rh)
	=\int\Psi^\vx(\uh)\,d\rh(\uh),\quad
\Psi^\edge(\bar\rh^\pi)
	=\int\Psi^\edge(\uh^\pi)\,d\rh(\uh)
\eeq
The Bethe ansatz is that $\log\Psi(\uh)$, for $\uh$ varying in (a possibly restricted subset of) $\splx^d$, is maximized at the \emph{replica symmetric} solution $h^j\equiv h$, with value $\log\Psi(h,\ldots,h)=\Phi(h)$. Condition \eqref{cond:ubd} says that this holds for all $\uh\in\supp\rh$ in an averaged sense. Of several natural modifications of $\rdc$ which we considered for the case of $d$ odd, all fail condition \eqref{cond:ubd}.

\section{Uniformly random $d$-regular graphs}\label{s:unif}

\subsection{Expectation of the partition function}

For $\usi$ a spin configuration on $G=(V,E)$, the edge empirical measure $L^\edge_n\in\splxpr$ of $\usi$ is defined by
$$L^\edge_n(\si,\si')
\equiv\f{1}{2\abs{E}} \sum_{(ij)\in E}
	[\kd{(\si_i,\si_j)}{(\si,\si')}
	+\kd{(\si_j,\si_i)}{(\si,\si')}
	],$$
and we write $L^\vx_n$ for its one-point marginal. Let $\splxpr^{d,n}$ denote the set of edge empirical measures associated to spin configurations on graphs in $\cG^\config_{d,n}$ (in the notation of \S\ref{sss:intro.unif.exp}).

\blem\label{l:unif.h}
For $\bh\in\splxpr^{d,n}$,
$$\abs{\spins}^n\P^\config_{d,n}[L^\edge_n=\bh]
= \exp\{n[ (d/2) H(\bh)-(d-1)H(\vh)]+\err\}$$
where $\abs{\err}\le c\log n$ for $c$ a finite constant depending only on $d,q$.

\bpf
Writing $\spins\equiv[q]$, we compute
\beq\label{e:p.h}
\abs{\spins}^n \P^\config_{d,n}[L^\edge_n=\bh]
=\CC(\bh)\MM(\bh)
\eeq
where
$$\CC(\bh)
\equiv
{n\choose n\vh_1, \ldots, n\vh_q}
\prod_\si
	{nd\vh_\si \choose nd\bbh{\si}{1},\ldots,nd\bbh{\si}{q}}$$
is the number of ways to assign spin values to the $nd$ half-edges subject to pair empirical measure $\bh$, and
$$\MM(\bh)
\equiv
[(nd-1)!!]^{-1}
	\prod_\si (nd \bbh{\si}{\si}-1)!!
	\prod_{\si\ne\si'} \sqrt{(nd\bbh{\si}{\si'})!}$$
is the number of perfect matchings $\gm\in\cM_{d,n}$ on half-edges respecting the spin assignment divided by the total number $\abs{\cM_{d,n}}$ of matchings of $[nd]$. With the convention $0^0=1$, Stirling's approximation (see e.g.\ \cite{\ww})
$$\Gam(z+1)=e^{O((z+1)^{-1})}[\sqrt{2\pi z}]^{\Ind{z>0}}(z/e)^z$$
gives
$\CC(\bh)
=\exp\{n  [dH(\bh)-(d-1)H(\vh)] + O(\log n)\}$.
Similarly, for $n$ even,
$$(n-1)!!
= \f{2^{n/2}}{\sqrt{\pi}} \Gam\Big(\f{n+1}{2}\Big)
= e^{O((n+1)^{-1})} (n/e)^{n/2} \sqrt{2}
$$
so $\gM(\bh)=\exp\{-(nd/2)H(\bh) + O(\log n)\}$, and the lemma follows.
\epf
\elem

\bpf[Proof of Thm.~\ref{t:fm.unif}]
The inequality follows trivially from Jensen's inequality. To compute $\E^\config_{d,n}[Z_n]$, let $\usi$ be drawn from the uniform distribution on $\spins^{[n]}$ and independently let $G_n\sim\P^\config_{d,n}$. Writing $L^\edge_n$ for the edge empirical measure of $\usi$ regarded as a spin configuration on $G_n$, the expected partition function for the model specified by $\upsi$ on $G_n$ can be expressed as
\begin{align}
\nonumber
\E^\config_{d,n}[Z_n]
&= \sum_{\bh\in\splxpr^{d,n}}
	\abs{\spins}^n
	\P^\config_{d,n}[L^\edge_n=\bh]
	\exp\{
	n\angl{\vxi}_\vh
	+(nd/2)\angl{\xi}_{\bh}\}\\
\label{e:unif.calc}
&=\sum_{\bh\in\splxpr^{d,n}}
	\exp\{n\bPhi(\bh) + O(\log n)\},
\end{align}
where the last line follows from \eqref{e:Phi.var} and Lem.~\ref{l:unif.h}. By Propn.~\ref{p:var}~\eqref{p:var.a}, $\bPhi$ attains its global maximum at an interior point $\bh^\star\in\splxpr$, which by \eqref{e:p.h} lies within distance $O(1/n)$ of $\splxpr^{d,n}$. Moreover the cardinality of $\splxpr^{d,n}$ is trivially bounded above by $(nd)^{q^2}$, so we find $\log \E^\config_{d,n}[Z_n]=n\bPhi(\bh^\star)+O(\log n)$, implying the theorem.
\epf

Turning to the proof of Cor.~\ref{c:fm.unif}, we first note that in the uniformly random $d$-regular graph ensemble, for any fixed $\ep>0$ and $t\ge0$ we have $\nt_t(G)\le\ep$ with overwhelming probability as $n\to\infty$:

\blem\label{l:unif.tree}
For any $\ep>0$, $t\ge0$ there exists a constant $\al\equiv \al(d,\ep,t)$ such that both $\P^\config_{d,n}[\nt_t(G)\ge\ep]$, $\P_{d,n}[\nt_t(G)\ge\ep]$ are bounded above by $e^{-\al\, n\log n}$.

\bpf
Consider the process of revealing the graph $G\sim\P^\config_{d,n}$ edge by edge. At the $k$-th step define $I_k$ to be the indicator that the edge forms a cycle of length $\le 2t$ within the graph revealed so far: if $\nt_t(G)\ge\ep$ then we must have $\sum_k I_k\ge \al_0 n$ for some $\al_0\equiv\al_0(d,\ep,t)>0$. For $k\le n(d-\al_0)/2$, the conditional probability of $I_k=1$ is at most $\al_1/n$ for $\al_1\equiv\al_1(d,\ep,t,\al_0)>0$. By a classical martingale inequality (see e.g.~\cite[Thm.~6.1]{\mcdiarmid}),
\begin{align*}
&\P^\config_{d,n}[\nt_t(G)\ge\ep]
\le\P^\config_{d,n}\Big[\sum_k I_k \ge \al_0 n\Big]
\le\P\Big[
	\sum_{k\le n(d-\al_0)/2} I_k \ge \f{\al_0 n}{2}
	\Big]\\
&\le\Big(\f{\al_1/n}{\al_0/(d-\al_0)}
		\Big)^{n\al_0/2}
	\Big(\f{1}{1-\al_0/(d-\al_0)}
		\Big)^{n(d/2-\al_0)}
\le e^{ - \al_2 \, n\log n}
\end{align*}
for $\al_2\equiv\al_2(d,\ep,t,\al_0,\al_1)>0$, which proves the result for $\P^\config_{d,n}$. The result for $\P_{d,n}$ follows immediately by \eqref{e:config.err}.
\epf
\elem

\bpf[Proof of Cor.~\ref{c:fm.unif}]
For $y\in\R$ let $p^\config_{d,n}(y)\equiv\P^\config_{d,n}[n^{-1}\log Z_n\ge\Phi+y]$. Take $\de>0$ such that $\liminf_{n\to\infty}[\Phi-\phi^\config_{d,n}]\ge2\de>0$, and recall from \S\ref{sss:intro.unif.exp} that the Azuma--Hoeffding bound implies $p^\config_{d,n}(-\de)\le e^{-c\de^2 n}$. For any $\ep>0$, it holds for sufficiently large $n$ that
\begin{align*}
e^{(\Phi-\ep)n}
&\le\E^\config_{d,n}[Z_n]
\le e^{(\Phi-\de)n}
+ \E^\config_{d,n}[Z_n \Ind{Z_n \ge e^{(\Phi-\de)n}}]\\
&\le e^{(\Phi-\de)n}
	+ p^\config_{d,n}(x) e^{Cn}
	+ (p^\config_{d,n}(-\de)-p^\config_{d,n}(x))
	e^{(\Phi+x)n},
\quad x\ge-\de.
\end{align*}
Taking $0<\ep < \de\wedge (c\de^2/2)$ and $0<x < (c\de^2/2)\wedge(C-\Phi)$ gives $p^\config_{d,n}(x)\ge e^{-n[C-\Phi+2\ep]}\equiv e^{-C_x n}$ which proves \eqref{e:unif.tree.above}. By Lem.~\ref{l:unif.tree} we can let $\ep\decto0$ and $t\incto\infty$ slowly enough in $n$ that
$$\lim_{n\to\infty}
	\f{\P^\config_{d,n}[\nt_t(G)\ge\ep]}
	{p^\config_{d,n}(x)}=0,$$
implying the existence of $G_n\in\cG^\config_{d,n}$, $G_n\lpc\treereg_d$ with $\liminf_{n\to\infty} n^{-1}\log Z_n\ge \Phi+x$.
\epf

\subsection{The Potts Bethe functional}

In the remainder of \S\ref{s:unif} we study global and local maxima of the Bethe functional $\bPhi$ for the ferromagetic Potts specification.

The Potts Bethe recursion for $\be,B\ge0$ preserves the subspace $\splxbal$ of measures
\beq\label{e:bias}
h\in\splx,\quad
h=(h_1,\ldots,h_q)
=\f1q ( 1+(q-1)b,1-b,\ldots,1-b),
\eeq
parametrized by $0\le b\le 1$. The map $\BP$ restricted to this subset is simply a univariate recursion $b\mapsto \bpb(b)$: in terms of the log-likelihood ratio $r\equiv\log (h_1/h_2)$, it has the particularly simple form
\beq\label{e:potts.bp.r}
\wt\bpb:r\mapsto B+(d-1) \log \f{e^{\be+r}+q-1}{e^r+e^\be+q-2},
\eeq
analyzed e.g.~in~\cite[Lem.~4.6]{\dms}. The maximal and minimal fixed points (in $b$) are given by
$$b^\free\equiv\lim_{t\to\infty}\bpb^{(t)}(0),
\quad b^\plus\equiv\lim_{t\to\infty} \bpb^{(t)}(1).$$
By monotonicity, $b^\free$ and $b^\plus$ are well-defined with $b^\free\le b^\plus$, and they are the only fixed points within this set.

In \S\ref{ss:potts.global.max} we use the result of Thm.~\ref{t:fm.unif} to prove Thm.~\ref{t:potts.sup} that $\bPhi$ attains its \emph{global} maximum on one of the edge empirical measures corresponding (via \eqref{e:bij}) to $b^\free,b^\plus$, and moreover that when $B>0$ these are the only possible global maximizers. This gives an essentially explicit solution to the Bethe variational problem for the ferromagnetic Potts model, and we do not know of a proof which does not go through the probabilistic results of Thm.~\ref{t:fm.unif}.

In \S\ref{ss:potts.loc.max} we supplement Thm.~\ref{t:potts.sup} with a study of the local maximizers of $\bPhi$. In view of the calculation \eqref{e:unif.calc}, we expect local maximizers of $\bh$ of $\bPhi$ to have the following probabilistic interpretation, which is in the spirit of results of~\cite{\mms}: if $G_n\lpc\treereg_d$, the factor model on $G_n$ conditioned to the subspace of configurations with edge empirical measure $L^\edge_n$ close to $\bh$ should converge locally weakly to the (Bethe) Gibbs measure corresponding to $\bh$. With this motivation in mind we classify the stationary points of $\bPhi$ and study which ones can be local maximizers.

\subsubsection{Global maximum}
\label{ss:potts.global.max}

We first review the well-known random-cluster ($\acr{fk}$) representation of the Potts model. The \emph{Edwards--Sokal} ($\acr{es}$) \emph{measure} on a finite graph $G=(V,E)$ is the probability measure on pairs $(\usi,\ueta)$, where $\usi\in\spins^V$ is a spin configuration as before and $\ueta\in\set{0,1}^E$ is a \emph{bond configuration}, given by
$$\es_G(\usi,\ueta)
\propto \prod_{i\in V} e^{B\kd{\si_i}{1}}
	\prod_{e=(ij)\in E}
	[ (1-p)^{1-\eta_e}+p^{\eta_e} \kd{\si_i}{\si_j} ],\quad
	p=1-e^{-\be}.$$
The marginal on $\usi$ is the Potts model with parameters $(\be,B)$, while the marginal on $\ueta$ is the \emph{$\acr{fk}$ measure}
$$\pi_G(\ueta)
\propto
\prod_{e\in E} p^{\eta_e}(1-p)^{1-\eta_e}
	\prod_{C\in\mathscr{C}(\ueta)}
	(1+(q-1)e^{-B\abs{C}}),$$
where the second product is taken over the collection $\mathscr{C}(\ueta)$ of connected components $C$ of $\ueta$ (with $\abs{C}$ the number of vertices in $C$). Conditioned on an $\acr{fk}$ configuration $\ueta$ with connected components $C_1,\ldots,C_k$ (with $k\equiv k(n)\le n$), a realization of $\usi$ from $\es_G(\cdot\giv\ueta)$ is obtained by giving the same spin $\si_\ell$ to all the vertices of each component $C_\ell$, independently over the different components, such that
$$\es_G(\si_\ell=\si\giv\ueta)
=u_\ell e^{B\abs{C_\ell}\kd{\si_\ell}{1}},\quad
u_\ell\equiv \f{1}{e^{B\abs{C_\ell}}+q-1}.$$

\bpf[Proof of Thm.~\ref{t:potts.sup}]
We assume without loss that $B>0$, with the result for $B=0$ following by continuity. Take $G_n$ any graph on $[n]$ and $\ueta$ any bond configuration on $G_n$, with connected components $C_1,\ldots,C_k$. For $\si\ne1$ let $Y_\ell\equiv Y^\si_\ell\equiv \abs{C_\ell}[\Ind{\si_\ell=\si}-u_\ell]$. It is easily verified that the cumulant generating functions $\ka_\ell(t)\equiv\log\E[e^{t Y_\ell}]=\log(1+u_\ell(e^{t\abs{C_\ell}}-1))-tu_\ell\abs{C_\ell}$ satisfy $\sup_{t\le B/2}\ka_\ell''(t)=\ka_\ell''(B/2)\le c$ for some finite constant $c\equiv c(B,q)$ \emph{not} depending on $\abs{C_\ell}$, and so
$$\vpi_{G_n}
\Big( \sum_{\ell=1}^k Y_\ell \ge \ep n \giv \ueta\Big)
	\le e^{-n[\ep t - c t^2/2]} \le e^{-n\ep^2/(2c)}
\quad\text{provided } \ep\le cB/2.$$
Thus it holds with $\es_{G_n}(\cdot\giv\ueta)$-probability at least $1-e^{-\sqrt{n}}$ that $L^\vx_n(\usi)$ belongs to the subspace $\splxbal_n$ of measures of $\splx$ within distance $n^{-1/8}$ of the space $\splxbal$ defined above. Consequently, if $Z^{\mathrm{bal}}_n$ denotes the Potts partition function of $G_n$ restricted to $\set{\usi:L^\vx_n(\usi)\in\splxbal_n}$, then $Z^{\mathrm{bal}}_n/Z_n\ge 1-e^{-\sqrt{n}}$ since this ratio is simply the average of $\es_{G_n}(L^\vx_n(\usi)\in\splxbal_n\giv\ueta)$ (as a function of $\ueta$) with respect to $\pi_{G_n}$.

Now recall the calculation \eqref{e:unif.calc} for the random regular graphs. For $\bh\in\splxpr$ with $\vh\notin\splxbal$, for sufficiently large $n$ we have $\vh\notin\splxbal_n$, so the contribution to $\E^\config_{d,n}[Z_n]$ from configurations $\usi$ with $L^\edge_n(\usi)=\bh$ is
$$\exp\{n\bPhi(\bh)+O(\log n)\}
\le e^{-\sqrt{n}}\E^\config_{d,n}[Z_n]
=\exp\{n\Phi-\sqrt{n}+O(\log n)\},$$
so we see that any global maximizer $\bh$ for $\bPhi$ must lie in $\splxbal$. Let $h\in\splxstar$ correspond to $\bh$ via \eqref{e:bij}: summing \eqref{e:bij} over $\si'\in[q]$ gives
$$\bzh\vh_\si=h_\si [(e^\be-1)h_\si +1],$$
which implies (since the right-hand side is increasing in $h_\si$ for $h_\si>0$) that $h$ is symmetric among the spins $\ne1$ and has a non-negative bias towards spin $1$. It is easily checked that the only such $h$ are $h^\free$ and $h^\plus$ (see e.g.~\cite[Lem.~4.6]{\dms}) which concludes the proof.
\epf

\subsubsection{Local maxima} \label{ss:potts.loc.max}

For the $q$-Potts model with $B\ge0$ and $\be>0$ we reparametrize $m\equiv e^B\ge1$, $\te\equiv1/(e^\be-1)>0$, so that \eqref{e:bpval} simplifies to
$$m \bpff(h_1)
=\bpff(h_2)=\ldots=\bpff(h_q)
	=\zh,\quad
\bpff(x)\equiv x^{-1}(x/\te+1)^{d-1}.$$
For $v\equiv\te/(d-2)$, $\bpff_-\equiv \bpff|_{(0,v]}$ is monotone decreasing while $\bpff_+\equiv \bpff|_{[v,1]}$ is monotone increasing, so clearly $\abs{\set{h_2,\ldots,h_q}}\le2$. More precisely, we have the following classification:

\bppn\label{p:potts.sols}
For the $q$-Potts model with parameters $m\equiv e^B\ge1$ and $\te\equiv1/(e^\be-1)>0$, for any $h\in\splxstar$ there exists $1\le\ell\le q$ and $\pi\in S_q$ such that
\beq\label{e:ltype}
h_{\pi(\si)}=\begin{cases}
	\QQ_\pm\equiv \bpff_\pm^{-1}(\zh/m), & \si=1,\\
	\pp_+\equiv \bpff_+^{-1}(\zh), & 2\le \si\le\ell,\\
	\pp_-\equiv \bpff_-^{-1}(\zh), & \ell+1\le \si\le q,
	\end{cases}
\eeq
and $\pi(1)=1$ if $m>1$.
\bnm[(a)]
\item \label{p:potts.onetype}
If $v\ge1$ or $m \bpff(v)>\bpff(1)$ then $\ell=1$ for all $h\in\splxstar$.
\item \label{p:potts.ltype} We say that $h\in\splxstar$ is an \emph{$\ell_\pm$-type solution} if \eqref{e:ltype} holds with $h_{\pi(1)}=\QQ_\pm$. For $\ell\ge2$, if $qv<1$ and $1\le m\le m_\ell(\te)$ then $\splxstar$ has $\ell_\pm$-type elements.
\enm

\bpf
It is clear from the preceding discussion that every $h\in\splxstar$ is of the form described in \eqref{e:ltype}.\footnote{The terminology degenerates in some cases, in particular when $m=1$: in this case $\QQ\in\set{\pp_\pm}$ so the $\ell_+$-type solutions coincide with the $(\ell+1)_-$-type solutions for $1\le\ell<q$, and the only $1_-$-type or $q_+$-type solution is the uniform distribution on $[q]$.} Fixed points of \eqref{e:potts.bp.r} correspond to $1_\pm$- or $q_\pm$-type solutions.

\medskip\noindent
\eqref{p:potts.onetype} If $v\ge1$ then $\bpff=\bpff_+$ is injective on $(0,1]$ so necessarily $\ell=1$. If $h\in\splxstar$ has $\ell>1$ then $\pp_+\le1/(\ell-1)$, so for $\QQ_\pm=\bpff_\pm^{-1}[\bpff(\pp_+)/m]$ to be well-defined we must have $m\bpff(v)\le F[1/(\ell-1)]$. In particular, if $m\bpff(v)>\bpff(1)$ then again all solutions must have $\ell=1$.

\medskip\noindent
\eqref{p:potts.ltype}
Assuming $v<1$, the function
\begin{align*}
g^{\ell,m}_\pm(p)
&\equiv \QQ_\pm(p) + (\ell-1)p+(q-\ell) \pp_-(p) \\
&\equiv \bpff_\pm^{-1}[\bpff(p)/m]+(\ell-1)p+(q-\ell)\bpff_-^{-1}[\bpff(p)]
\end{align*}
is well-defined for $p\in[p_0,1]$ where $p_0\equiv p_0(m)\equiv \bpff_+^{-1}[m\bpff(v)]$. Note that $g^{\ell,m}_\pm(1)>\ell-1$, and $\lim_{m\decto 1}p_0(m)=v$ which implies $\lim_{m\decto1} g^{\ell,m}_\pm(p_0) = qv$. If $\ell\ge2$ and $qv<1$ then continuity of $g^{\ell,m}_\pm(p)$ implies that for $m\ge1$ sufficiently small we will have $g^{\ell,m}_\pm(\pp_+)=1$ for some $\pp_+\in[p_0,1]$, giving an $\ell_\pm$-type solution as claimed.
\epf
\eppn

We next study which of the stationary points classified in Propn.~\ref{p:potts.sols} correspond to local maximizers for $\bPhi$.

\bppn\label{p:potts.locmax}
In the setting of Propn.~\ref{p:potts.sols},
\bnm[(a)]
\item \label{p:potts.notlocmax}
Solutions of type $\ell_\pm$ with $\ell>2$ are never local maximizers.
\item \label{p:potts.hightemp.locmax} For $m\ge1,\te>0$ both sufficiently small, there exist both $1_+$-type and $2_-$-type solutions which are strict local maximizers with (strictly) negative-definite Hessians.\footnote{If $qv<1$ then there can be no $1_-$-type solutions.}
\enm

\bpf
\medskip\noindent
\eqref{p:potts.notlocmax}
Let $\bh\in\splxpr$ be the stationary point of $\bPhi$ corresponding to $h$ via \eqref{e:bij}. We will apply the correlation criterion \eqref{e:corr} with $\bbph{\si}{\si'}\equiv \vph_\si +\vph_{\si'}$. Let $\bh$ correspond to $h\in\splxstar$ via \eqref{e:bij}, so that
$$\bh_{\si'\giv \si}
\equiv \f{\bbh{\si}{\si'}}{\vh_\si }
= \f{h_{\si'}(\te+\kd{\si}{\si'})}{\te+h_\si }.$$
If we assume $\angl{\vph}_h=0$, then
$$\E[\vph_Y \giv X=\si]
= \f{\te}{\te+h_\si } \angl{\vph}_h
	+ \f{h_\si }{\te+h_\si } \vph_\si 
= \gam_\si \vph_\si $$
for $\gam_\si\equiv h_\si /(\te+h_\si )$. Thus $\E[\bbph{X}{Y} \giv X]=(1+\gam_X)\vph_X$, and \eqref{e:corr} becomes
\begin{align}
\nonumber &2(\E\vph_X)^2
\ge  \E\Big[(1+\gam_X)\vph_X^2
	\f{(d-1)\gam_X-1}{d-2}
	\Big]\\
\label{e:corr.simp}
&= \E\Big[
	(1+\gam_X)\vph_X^2 \f{h_X-v}{\te+h_X}\Big]
=\f{1}{\te+\nrm{h}^2}
	\sum_\si h_\si (1+\gam_\si)(h_\si -v)\vph_\si ^2
\end{align}
(using $\vh_\si =h_\si (\te+h_\si )/(\te+\nrm{h}^2)$ for the last identity). If $\bh$ is an $\ell$-type solution with $\ell>2$ then $\vph_\si =\kd{\si}{\pi(2)}-\kd{\si}{\pi(3)}$ (for $\pi$ as in \eqref{e:ltype}) clearly violates \eqref{e:corr.simp}, so $\bh$ cannot be a local maximizer of $\bPhi$.

\medskip\noindent
\eqref{p:potts.hightemp.locmax}
Let $m=1$ and $\te$ sufficiently small so that a $1_+$-type (and $2_-$-type) solution $h\equiv (\QQ_+,\pp_-,\ldots,\pp_-)\in\splxstar$ exists, given by taking the log-likelihood ratio $\rr\equiv\log(\QQ_+/\pp_-)$ to be the maximal fixed point of the mapping $\wt\bpb$ of \eqref{e:potts.bp.r}. For $d\ge3$ and $0<\ep\le1$,
$$\wt\bpb[(d-1-\ep)\be]
\ge(d-1) (\be-\log q)
>
(d-1-\ep)\be\quad\forall
\be
>
\f{(d-1)\log q}{\ep},$$
so crudely we have $\rr\ge(3/2)\be$ for all $\be\ge2(d-1)\log q$. Let $\bh\in\splxpr$ be the stationary point corresponding to this fixed point: recalling \eqref{e:locmax.delta}, for $\bde\in\splxprpm$ we calculate
\begin{align*}
\f{1}{\bzh} \angl{(\vde/\vh)^2}_{\vh}
&= \f{\vde_1^2}
		{\QQ_+(e^\be \QQ_+ + (q-1) \pp_-)}
	+
	\f{\sum_{\si\ne1}\vde_\si^2}
		{\pp_-(\QQ_++(e^\be+q-2)\pp_-)}
\le \f{\vde_1^2}{e^\be \QQ_+^2}
	+ \f{\sum_{\si\ne1}\vde_\si^2}{\QQ_+\pp_-},\\
\f{1}{\bzh} \angl{(\bde/\bh)^2}_{\bh}
&\ge
	\f{\bbde{1}{1}^2}{\QQ_+^2 e^\be}
	+\f{2\sum_{\si\ne1}\bbde{1}{\si}^2}{\QQ_+\pp_-}
	+\f{ \sum_{\si,\si'\ne1} \bbde{\si}{\si'}^2
		}{e^\be \pp_-^2}.
\end{align*}
Since $\pp_-\le e^{-(3/2)\be}\QQ_+$ for sufficiently large $\be$,
$$\lim_{\be\to\infty}
	\f{e^\be\pp_-^2}{\bzh} \angl{(\vde/\vh)^2}^2_{\vh} =0,\quad
\liminf_{\be\to\infty}
	\f{e^\be\pp_-^2}{\bzh} \angl{(\bde/\bh)^2}_{\bh}
	\ge\sum_{\si,\si'\ne1}\bbde{\si}{\si'}^2,$$
so for any fixed $\ep>0$ we have $\pd_\eta^2\Phi_\mu(\bh+\eta\bde)|_{\eta=0}<0$ uniformly over all $\bde\in\splxprpm$ with $(q-1)^2\sum_{\si,\si'\ne1}\bbde{\si}{\si'}^2\ge\ep^2$ once $\be$ is sufficiently large (depending on $\ep$).

Suppose instead $(q-1)^2\sum_{\si,\si'\ne1}\bbde{\si}{\si'}^2\le\ep^2$: by Cauchy--Schwarz $\sum_{\si,\si'\ne1}\abs{\bbde{\si}{\si'}}\le\ep$, so
\begin{align*}
&\limsup_{\be\to\infty}
\f{\QQ_+\pp_-}{\bzh}
\Big[2(d-1)\angl{(\vde/\vh)^2}_{\vh}
-d\angl{(\bde/\bh)^2}_{\bh}\Big]
\le 2(d-1)\sum_{\si\ne1}\vde_\si^2
	-2d\sum_{\si\ne1}\bbde{1}{\si}^2\\
&\le 2(d-1)\sum_{\si\ne1} [\abs{\bbde{1}{\si}}+\ep]^2
	-2d\sum_{\si\ne1} \bbde{1}{\si}^2
\le -2\sum_{\si\ne1} \bbde{1}{\si}^2
	+ 2(d-1)
	\Big[
	2\ep\sum_{\si\ne1}\abs{\bbde{1}{\si}} + (q-1)\ep^2
	\Big].
\end{align*}
On the other hand, $\bde\in\splxprpm$ implies
$$
2\Big|\sum_{\si\ne1} \bbde{1}{\si}\Big|
=\abs{\bbde{1}{1}+\ep}
\ge \abs{\bbde{1}{1}}-\ep
\ge \Big[
	1-2\sum_{\si\ne1} \bbde{1}{\si}^2
	-\f{\ep^2}{(q-1)^2}
	\Big]^{1/2}-\ep$$
so by choosing $\ep>0$ sufficiently small we can guarantee that for $\be$ large enough, $\pd_\eta^2\Phi_\mu(\bh+\eta\bde)|_{\eta=0}<0$ uniformly over all $\bde\in\splxprpm$, implying that $\bh$ is a strict local maximizer of $\bPhi$ with strictly negative-definite Hessian.

This concludes the proof for $m=1$, and the conclusion for $m>1$ sufficiently small follows by a perturbative argument: arguing similarly as in the proof of Propn.~\ref{p:potts.sols}~\eqref{p:potts.ltype}, for $1\le m<m_0$ the equations
\begin{align*}
g^{1+}(p)&\equiv \bpff_+^{-1}[\bpff(p)/m]+(q-1)p=1,\\
g^{2-}(p)&\equiv \bpff_-^{-1}[\bpff(p)/m]+\bpff_+^{-1}[\bpff(p)]+(q-2)p=1
\end{align*}
have solutions $\pp_-^{1+}(m),\pp_-^{2-}(m)$, corresponding to $1_+$-type and $2_-$-type solutions respectively, which are continuous in $m$ with initial values $\pp^{1+}(1)=\pp^{2-}(1)=\pp_-$ corresponding to the solution considered above at $m=1$. For sufficiently small $m$, it follows by continuity that the Hessians at the stationary points $\bh^{1+}(m),\bh^{2-}(m)$ corresponding to $\pp^{1+}_-(m),\pp^{2-}_-(m)$ will be strictly negative-definite, implying strict local maximizers as claimed.
\epf
\eppn

\brmk
Related to the study of local maxima is the question of the local stability of the Bethe recursion. For the Potts specification \eqref{e:potts}, the linear (differential) mapping $D_h\equiv D\BP(h)$ defined on the space $\set{\de:\sum_\si\de_\si=0}$ by
$$D_h\de\equiv\lim_{\eta\to0}
	\f{\BP(h+\eta\de)-\BP(h)}{\eta}$$
can be explicitly diagonalized when $h\in\splxstar$ and shown to have all eigenvalues positive, with maximal eigenvalue greater than $1$ at $\ell_\pm$-type solutions with $\ell>2$ and at $2_+$-type solutions. At a $2_-$-type solution (assuming $m>1$, so it is not also a $1_+$-type solution) the maximal eigenvalue is less than $1$ if and only if
\beq\label{e:stab}
\f{\pp_+^2}{\pp_+-v}
> \f{d-2}{d-1} + \f{\QQ_-^2}{v-\QQ_-}
+(q-2)\f{\pp_-^2}{v-\pp_-}.
\eeq
However, if $h$ is not the uniform measure on $[q]$ then $D_h$ is \emph{not} symmetric and so does not have orthonormal eigenbasis, so having all eigenvalues less than $1$ need not imply contractivity of $D_h$. It is not clear how to relate \eqref{e:stab} to the local stability of the non-linear map $\BP$.
\ermk

\section{Recursive graph decomposition}\label{s:decomp}

In this section we prove Thm.~\ref{t:fm}. Recall from \S\ref{ss:intro.gen} the notation $\oo(t,x)$; we also let $c$ denote a finite positive constant which is permitted to change from line to line but depends only on $d,q,\upsi$. The following lemma, whose proof we defer to the end of the section, reduces the free energy computation on general $G_n\lpc\treereg_d$ to the case of $d$-regular graphs.

\blem\label{l:dreg}
If $G_n\lpc\treereg_d$ with $d$ even, then there exists a $d$-regular (multi-)graph sequence $G_n'\lpc\treereg_d$ with free energy $\phi_n'$ such that $\lim_{n\to\infty}(\phi_n-\phi_n')=0$.
\elem

We prove Thm.~\ref{t:fm} via the following propositions about the operation $\rdc$ on $d$-regular graphs.

\bppn\label{p:reduc.ubd}
Let $G=(V,E)$ be any finite $d$-regular graph, and recall $\E_G$ denotes expectation over the uniformly random vertex $I_G$. Under the conditions of Thm.~\ref{t:fm},
$$\E_G[\log Z_G-\log Z_{\rdc G}]
\le\Phi+\oo(t,\nt_t(G)).$$
\eppn

\bppn\label{p:unif.tree}
Suppose $G_n\lpc\treereg_d$. Then for all $\ep_0>0$ and all $t\ge0$,
$$\lim_{n\to\infty}
	\E_n\Big[
	\max_{0\le j\le (1-\ep_0)n}
	\nt_t(\rdc^j G_n) \Big]
	= 0.$$
\eppn

We first assume the preceding results and derive Thm.~\ref{t:fm}:

\bpf[Proof of Thm.~\ref{t:fm}]
Take $\ep_0>0$ fixed, let $n_0\equiv\flr{(1-\ep_0)n}$, and let $\ep\equiv\ep_n$ be defined by $n_0\equiv(1-\ep_n)n$. Express the free energy of the factor model on $G_n$ as the telescoping sum
$$\phi_n
= \f1n\sum_{j=0}^{n_0-1}
	\E_n[\log Z_{\rdc^j G_n}-\log Z_{\rdc^{j+1}G_n}]
	+ \f1n\E_n[\log Z_{\rdc^{n_0}G_n}].$$
By definition of $\rdc$, $\rdc^{n_0}G_n$ is a $d$-regular graph on $\ep n$ vertices, so $n^{-1}\abs{\log Z_{\rdc^{n_0}G_n}}\le c\ep$. Next, Propn.~\ref{p:reduc.ubd} gives
$$\max_{0\le j<n_0}
\E_n[\log Z_{\rdc^j G_n}-\log Z_{\rdc^{j+1}G_n}]
	-\Phi
\le \oo\big(t,
	\max_{0\le\ell<n_0}
	 \E_n[\nt_t(\rdc^\ell G_n)] \big)$$
(where we may freely move the expectation inside $\oo(t,\cdot)$ by uniform boundedness of $\oo$). In the limit $n\to\infty$ the right-hand side above tends to $\oo(t)$ by Propn.~\ref{p:unif.tree}, so the telescoping sum yields
$$\limsup_{n\to\infty}\phi_n
\le (1-\ep_0)[\Phi+\oo(t)] + c\ep_0.$$
The result follows by taking first $t\to\infty$ and then $\ep_0\decto0$.
\epf

The remainder of this section is devoted to proving Propns.~\ref{p:reduc.ubd} and~\ref{p:unif.tree} and Lem.~\ref{l:dreg}. Recall \eqref{e:perm.bd} that $\psimin^{d+1} \le Z_G/Z_{G_-}\le \psimax^{d+1}$; similarly
\beq\label{e:perm.bd.e}
\psimax^{d/2}\ge\f{Z_{G^\pi}}{Z_{G^-}}
\ge \psi(\si^\perm,\si^\perm)^{d/2} \numin(\usi_{\pd I}\equiv\si^\perm)
\ge \f{\psimin^{d/2}}{\abs{\spins}^d \ratio^{\abs{\pd I}} \ratio^{\abs{\pd(\pd I)}}}
\ge \f{\psimin^{d/2}}{\abs{\spins}^d \ratio^{d(d+1)}}
\eeq
for $\ratio\equiv\psimax/\psimin$.

\bpf[Proof of Propn.~\ref{p:reduc.ubd}]
Fix $t\ge0$, and let $\rh\equiv\rh(G,I,t)$ as in the statement of Thm.~\ref{t:fm}. It follows from condition \eqref{cond:prod} --- recalling \eqref{e:ratio.gen} and \eqref{e:ratio.prod} and making use of the boundedness of $\upsi$ --- that
$$\absb{\f{Z_G}{Z_{G_-}}
	-\Psi^\vx(\bar\rh)
	}
+\absb{\f{Z_{G^\pi}}{Z_{G_-}}
	-\Psi^\edge(\bar\rh^\pi)
	}
\le\oo(t,\nt_t(G))\quad
\text{with probability}
\ge 1-\oo(t,\nt_t(G)).$$
Then, by the choice of permutation in the definition of $\rdc$ and using the bounds \eqref{e:perm.bd}, \eqref{e:perm.bd.e} on the ratios $Z_G/Z_{G_-}$ and $Z_{G^\pi}/Z_{G_-}$, we find
\begin{align*}
\E_G
	\Big[\log \f{Z_G}{Z_{\rdc G}}\Big]
&=\E_G\Big[
	\log\min_\pi \f{Z_G/Z_{G_-}}{Z_{G^\pi}/Z_{G_-}}\Big]
\le\E_G\Big[
	\log
	\f
	{\Psi^\vx(\bar\rh)}
	{\Psi^{\edge,\sym}(\bar\rh)}
	\Big]
	+\oo(t,\nt_t(G))\\
&\le\E_G\Big[
	\log \max_{\uh\in\supp\rh}\Psi^\sym(\uh)
	\Big]
	+\oo(t,\nt_t(G)).
\end{align*}
Condition \eqref{cond:ubd} of Thm.~\ref{t:fm} gives that the first term is $\le\Phi+\oo(t)$ which implies the result.
\epf

\bpf[Proof of Propn.~\ref{p:unif.tree}]
Let $G$ be a graph on $n$ vertices. For fixed $\ep_1>0$ (to be determined), let $\ep\equiv\ep_n$ be defined by $\ep_n n \equiv \ceil{\ep_1 n}$, and consider $\ep n$ successive applications of $\rdc$. Assume $t\ge0$ is a (large) power of $2$, and for $v\in G$ let $\Gam_v$ be the indicator of the event that $v$ is one of the $n[1-\nt_t(G)]$ vertices with $\ball{t}{v}\cong\treereg^t_d$, but $\ball{t/2}{v}\not\cong\treereg^{t/2}_d$ in $\rdc^j G$ for some $j\le\ep n$. Then $\E[\Gam_v]$ is bounded above by the probability that at least $t/2$ vertices are deleted along a length-$t$ geodesic path started from $v$. For a single path, this probability is $\P[H\ge t/2]$ where $H$ is a hypergeometric random variable with parameters $n,t,\ep n$. By a standard hypergeometric tail bound (see~\cite{\chvatal}),
$$
\P[H\ge\gam\ep t]
\le \Big[ (1/\gam)^{\gam t/n}
	\Big(\f{1-t/n}{1-\gam t/n}\Big)^{1-\gam t/n}
	\Big]^{\ep n}
\le (1/\gam)^{\gam\ep t}
	e^{(\gam-1)\ep t}
\le \exp\{-\gam\ep t[\log\gam-1] \}.$$
Taking $\gam=1/(2\ep)$ and summing over $d^t$ geodesics gives
$$\E[\Gam_v]
\le \exp\Big\{
	-\f{t}{2} \Big( \log\f{1}{2\ep}-2\log d-1 \Big)
	\Big\},$$
which can clearly be made $\le e^{-t/4}$ by taking $\ep_1$ sufficiently small (depending only on $d$) and $n,t$ large. Markov's inequality applied to the $\Gam_v$ ($v\in G$) then gives
$$\P\Big[
\max_{0\le j\le\ep n}
	(n-j) \nt_{t/2}(\rdc^j G)
\ge n\nt_t(G) + ne^{-t/8}\Big]
\le e^{-t/8}.$$
Iterating for $L=\ceil{(\log\ep_0)/\log(1-\ep_1)}$ steps, and writing $t(L)\equiv t/2^L$, gives that with probability at least $1-L e^{-t(L)/8}$,
$$(n-j)\nt_{t/2^L}(\rdc^j G) \le n \nt_t(G) + L n e^{-t(L)/8}\quad
\text{for all }
0\le j\le (1-\ep_0)n.$$
Therefore, for $G=G_n\lpc\treereg_d$,
$$\E_n\big[
\max_{0\le j\le (1-\ep_0) n}
\nt_{t/2^L}(\rdc^j G_n)
\big]
\le
\f{n\E_n[\nt_t(G_n)] + 2Ln e^{-t(L)/8}}{\ep_0 n}.$$
In the limit $n\to\infty$ the right-hand side tends to $Le^{-t(L)/8}/\ep_0$ which decreases to zero as $t\to\infty$, but the left-hand side is non-decreasing in $t$ so it must in fact tend to zero as $n\to\infty$ for all $t$, as stated.
\epf

\bpf[Proof of Lem.~\ref{l:dreg}]
Fix $n$ for the moment and suppress it from the notation. Delete edges in $G$ incident to vertices of degree larger than $d$ until none remain, and denote the resulting graph $G''=(V,E'')$. Let $U$ denote the set of vertices incident to any edge in $E\setminus E''$; arguing as for the bounds \eqref{e:perm.bd} and \eqref{e:perm.bd.e} then gives
$$\psimax^{\abs{E\setminus E''}}
\ge \f{Z_G}{Z_{G''}}
\ge \psimin^{\abs{E\setminus E''}}
	\nu_{G''}(\usi_U\equiv\si^\perm)
\ge \f{\psimin^{\abs{E\setminus E''}}}
	{\abs{\spins}^{\abs{U}}
		\ratio^{\abs{U}+\abs{\pd U}}}
\ge \f{\psimin^{\abs{E\setminus E''}}}
	{\abs{\spins}^{\abs{U}}
		\ratio^{(d+1)\abs{U}}}$$
(where the last inequality uses that $G''$ has maximum degree at most $d$). Then note that $\abs{U}/2
\le\abs{E\setminus E''}
\le\sum_{v\in G} [(D_v-d)\vee0]
\le \abs{V}\,\E_G[D_{I_G} \Ind{D_{I_G}>d}]$,
so $\abs{\log Z_G-\log Z_{G''}}\le c\abs{V}\,\E_G[D_{I_G} \Ind{D_{I_G}>d}]$.

Now let $W$ denote the set of vertices $v\in V$ whose degree $D''_v$ in $G''$ is less than $d$, and add $d-D''_v$ new half-edges leaving from each such $v$. The number of unmatched half-edges resulting from this operation has the same parity as
$\sum_{v\in G} (D_v-d) = 2\abs{E}-d\abs{V}$,
so it is even. Taking a random matching of these half-edges results in a $d$-regular (multi-)graph $G'$. The total number $\de_E$ of edge insertion or deletion operations to go from $G''$ to $G'$ is at most a constant times $\abs{V}\E_G[(D_{I_G}\vee1)\Ind{D_{I_G}\ne d}]$, since if $v\in W$ then either $v\in U$ or $D_v<d$ in $G$. Each such operation changes the log-partition function by at most an additive constant, since all graphs involved have maximum degree at most $d$. Therefore
$$\abs{\log Z_{G''}-\log Z_{G'}}
\le c\abs{V}\,\E_G[(D_{I_G}\vee1)\Ind{D_{I_G}\ne d}].$$
Combining with the previous bound on $\abs{\log Z_G-\log Z_{G''}}$ gives
\begin{align*}
\abs{\phi_n-\phi_n'}
&\le c\,\E_n[(D_{I_n}\vee1)\Ind{D_{I_n}\ne d}]\\
&\le c L\, \P_n[D_{I_n}\ne d]
	+c\,\E_n[D_{I_n}\Ind{D_{I_n}\ge L}].
\end{align*}
The lemma follows by taking first $n\to\infty$ and then $L\to\infty$ in the bound above.
\epf

\section{The Potts free energy density}\label{s:potts}

In this section we prove our main result Thm.~\ref{t:potts} giving the free energy density of the $q$-Potts model on graphs converging locally to the $d$-regular tree with $d$ even. Let $\splxfk{\ep}$ denote the measures $h\in\splx$ of form \eqref{e:bias} with $b\in[b^\free-\ep,b^\plus+\ep]$.

\bppn\label{p:potts.opt}
For the Potts model with $\be,B\ge0$, $\log\sup_{\uh\in(\splxfk{0})^d} \Psi^\sym(\uh)=\Phi$.
\eppn

\bppn\label{p:potts.prod}
For the Potts model with $\be\ge0$ and $B>0$, for all finite $d$-regular graphs $G$, all $v\in G$, and all $t\in\Z_{\ge0}$, there is a measure $\rh\equiv\rh(G,v,t)$ on $(\splxfk{\oo(t)})^d$ such that for $I=I_G$ and $\numin\equiv\nu_{G\setminus I}(\usi_{\pd I}=\cdot)$, the measure $\bar\rh$ of\eqref{e:bar.rho} satisfies
$\E_G[\tv{\numin-\bar\rho}]\le\oo(t,\nt_t(G))$.
\eppn

\bpf[Proof of Thm.~\ref{t:potts}]
The result follows from Thm.~\ref{t:fm} since condition \eqref{cond:prod} is verified by Propn.~\ref{p:potts.prod} while condition \eqref{cond:ubd} is verified by Propn.~\ref{p:potts.opt} (using continuity of $\Psi^\sym$ on $\splx^d$).
\epf

\bpf[Proof of Propn.~\ref{p:potts.opt}]
{\it Step 1.} For $\uh\in(\splxfk{0})^d$ we abuse notation and write $\Psi(\uh)\equiv\Psi(\ub)$ with $h^j(1)\equiv(1+(q-1)b^j)/q$. Then
\begin{align*}
C^{-d} \Psi^\vx(\ub)
&= e^B \prod_{j=1}^d (1+\gam b^j)
+ (q-1)\prod_{j=1}^d(1-\gam b^j/(q-1)),\\
C^{-d/2} \Psi^\edge(\ub)
&= \prod_{j=1}^{d/2}(1+\gam b^{2j-1}b^{2j})
\end{align*}
where $C\equiv(e^\be+q-1)/q$, $\gam\equiv(q-1)(e^\be-1)/(e^\be+q-1)>0$. Both $\Psi^\vx$ and $\Psi^{\edge,\sym}$ are affine in each $b^j$ (keeping $(b^k)_{k\ne j}$ fixed), and so $\Psi^\sym$ is maximized with $b^j\in\set{b^\free,b^\plus}$. Therefore
$$\sup_{\uh\in(\splxfk{0})^d} \Psi^\sym(\uh)
=\sup_{\ub\in\set{b^\free,b^\plus}^d}\Psi^\sym(\ub).$$

\medskip \noindent {\it Step 2.} Let $g^\forp(b)\equiv\log\Psi(b,b^\forp,\ldots,b^\forp)-(d/2)\log C$ for $\forp\in\set{\free,\plus}$; we claim that $g^\forp$ is constant in $b$. To see this, note that $\Psi^\vx(\uh)=Z_S(\uh)$, the partition function of the Potts model on the star graph $S\equiv\treereg^1_d$ with boundary conditions $\si_j\sim h^j$ independently for the vertices $j\in\pd\rt$. Similarly, $\Psi^\edge(\uh)=Z_R(\uh)$, the partition function of the Potts model on the graph $R$ of $d/2$ disjoint edges $(2k-1,2k)$, again with boundary conditions $\si_j\sim h^j$ independently for all $j$. Now if $S'$ is $S$ with the edge $(\rt,1)$ disconnected, then, using the BP relations for $b^\forp$,
$$\f{Z_S(b,b^\forp,\ldots,b^\forp)}
{Z_{S'}(b,b^\forp,\ldots,b^\forp)}
=1+\gam b b^\forp.$$
But $Z_{S'}$ does not depend on $b$, so with $\ub^\forp\equiv(b^\forp,\ldots,b^\forp)$ we have
\beq\label{e:const.in.b}
\f{Z_S(b,b^\forp,\ldots,b^\forp)}
{Z_S(\ub^\forp)}
=\f{1+\gam b b^\forp}{1+\gam (b^\forp)^2}.
\eeq
If $R'$ denotes $R$ with the edge $(1,2)$ disconnected then by the same argument \eqref{e:const.in.b} holds with $S,S'$ replaced by $R,R'$, and so
$$\exp\{g^\forp(b)-g^\forp(b^\forp)\}
= \f{Z_S(b,b^\forp,\ldots,b^\forp)}{Z_R(b,b^\forp,\ldots,b^\forp)} \f{Z_R(\ub^\forp)}{Z_S(\ub^\forp)}=1.$$

\medskip \noindent {\it Step 3.} Let $\ub\in\set{b^\free,b^\plus}^d$, and let $\ell$ denote the number of indices $j$ for which $b^j=b^\plus$. Then
\begin{align*}
C^{-d}\Psi^\vx(\ub)
&= e^B(1+\gam b^\free)^d
	\Big(
	\f{1+\gam b^\plus}{1+\gam b^\free} \Big)^\ell
+(q-1)(1-\gam b^\free/(q-1))^d
	\Big( \f{1-\gam b^\plus/(q-1)}{1-\gam b^\free/(q-1)} \Big)^\ell\\
&\equiv A_0 e^{a_0\ell} +A_1 e^{-a_1\ell}
\equiv \exp\{f^\vx(\ell)\},\quad A_j,a_j>0.
\end{align*}
By Jensen's inequality,
$$\log \Psi^{\edge,\sym}(\ub)
\ge \f{1}{d!}\sum_{\pi\in S_d}
	\log \Psi^\edge(\ub^\pi)
\equiv f^\edge(\ell) + (d/2)\log C.$$
Note that for $\ell\in\set{0,1,d-1,d}$, both sides are equal to $\log\Psi^\edge(\ub)$. For $0\le\ell\le d$,
\begin{align*}
f^\edge(\ell) &= [2(d-1)]^{-1}
\{
\ell(\ell-1) C_{\plus\plus}+(d-\ell)(d-\ell-1) C_{\free\free}
+2\ell(d-\ell) C_{\plus\free}]
\}\\
&\equiv a_4\ell^2+a_3\ell+a_2,
\end{align*}
where $C_{\forp\forp'}\equiv\log(1+\gam b^\forp b^{\forp'})$ for $\forp,\forp'\in\set{\free,\plus}$, and $2(d-1)a_4=C_{\plus\plus}+C_{\free\free}-2C_{\plus\free}$ which for $b^\free<b^\plus$ is strictly positive by the arithmetic-geometric mean inequality. If we now consider $f\equiv f^\vx-f^\edge$ as a function of $\ell\in\R$, then
$$f'(\ell)
= \f{A_0 a_0 e^{a_0 \ell}-A_1 a_1e^{-a_1\ell}}
{A_0e^{a_0 \ell}+A_1 e^{-a_1\ell}}-2 a_4\ell-a_3$$
tends to $\mp\infty$ as $\ell\to\pm\infty$. Moreover
$$f'''(\ell)
= -
\f{A_0A_1 (a_0+a_1)^3 e^{(a_0-a_1)\ell}}{(A_0e^{a_0\ell}+A_1e^{-a_1\ell})^3}
(A_0e^{a_0\ell}-A_1e^{-a_1\ell})$$
which can have at most one real zero. Thus $f'$ has at most one inflection point, hence at most three real zeroes; further, if there are three zeroes then the middle one corresponds to a local minimum of $f$. But $f(\ell)=g^\forp(\ell)$ for $\ell\in\set{0,1,d-1,d}$, so Step 2 implies $f(0)=f(1)$ and $f(d-1)=f(d)$, and consequently $f'$ has zeroes in $(0,1)$ and $(d-1,d)$. Therefore $f$ cannot have a local maximum in $[1,d-1]$, so it is maximized over $\set{0,\ldots,d}$ with $\ell\in\set{0,d}$, completing the proof.
\epf

To decompose the measures $\numin$ as mixtures over $(\splxfk{\oo(t)})^d$, we use the random-cluster ($\acr{fk}$) representation reviewed in \S\ref{s:unif}.

\bpf[Proof of Propn.~\ref{p:potts.prod}]
For $0\le s\le t<\infty$ and $v\in G$ let $A_{s,t}(v)\equiv\ball{t}{v}\setminus\ball{s}{v}$. It holds with probability at least $1-\oo(t,\nt_t(G))$ that for $I=I_G$ with $\pd I=(v_1,\ldots,v_d)$, the balls $\ball{t}{v_i}$ defined with respect to $G_-\equiv G\setminus I$ are pairwise disjoint and isomorphic to the first $t$ levels of the $(d-1)$-ary tree. On this event, $\numin\equiv\nu_{G_-}(\usi_{\pd I}=\cdot)$ has the decomposition
$$\numin(\cdot)
=\sum_{\ueta_{\bA_{s,t}}}
	\esmin(\cdot\giv\ueta_{\bA_{s,t}})
	\pimin(\ueta_{\bA_{s,t}}),\quad
\bA_{s,t}\equiv\bigcup_{i=1}^d A_{s,t}(v_i)$$
(where $\esmin\equiv\es_{G_-}(\usi_{\pd I}=\cdot)$ and $\pimin\equiv\pi_{G_-}$). We claim that with $s=\log t$ the measures $\esmin(\cdot\giv\ueta_{\bA_{s,t}})\in\splx_{\spins^d}$ are approximately in $(\splxfk{\oo(t)})^d$. Indeed, write $\bB_t$ for the disjoint union of the balls $\ball{t}{v_i}$, and decompose
\beq\label{e:potts.fk}
\esmin(\cdot\giv\ueta_{\bA_{s,t}})
=\sum_{\ueta_{\bB_t}}
	\esmin(\cdot\giv\ueta_{\bB_t})
	\pimin(\ueta_{\bB_t}\giv\ueta_{\bA_{s,t}})
\eeq
For $\ddagger\in\set{0,1}$ write $\esmin^\ddagger$ and $\pimin^\ddagger$ for the $\acr{es}$ and $\acr{fk}$ measures respectively on $G_-$ conditioned on $\ueta_{G_-\setminus B_t}\equiv\ddagger$ (the measures $\esmin^0$, $\pimin^0$ restrict simply to the $\acr{es}$, $\acr{fk}$ on $\bB_t$). The $\acr{fkg}$ property of random-cluster measures (see e.g.\ \cite[Thm.~III.1(i)]{\bbck}) implies the stochastic domination relations
$
\pimin^0(\ueta_{\bB_t}=\cdot\giv\ueta_{\bA_{s,t}})
\preccurlyeq
\pimin(\ueta_{\bB_t}=\cdot\giv\ueta_{\bA_{s,t}})
\preccurlyeq
\pimin^1(\ueta_{\bB_t}=\cdot\giv\ueta_{\bA_{s,t}})$.
Now note that
$$\f{\pimin^0(\ueta_{\bB_t}
	\giv\ueta_{\bA_{s,t}})}
{\pimin^1(\ueta_{\bB_t}
	\giv\ueta_{\bA_{s,t}})}
\propto
\rh(\ueta_{\bB_t})\equiv \f{\prod_{C\in\mathscr{S}} [1+(q-1)e^{-B\abs{C}}]}
{1+(q-1) e^{-B\abs{C_\infty}}}$$
where $\mathscr{S}$ is the set of connected components $C$ of $\ueta_{\bB_t}$ joining $\ball{s}{v_i}$ to the boundary of $\ball{t}{v_i}$ for some $i$, and $C_\infty$ is the union of these components with the complement of $\bB_t$. Each $C\in\mathscr{S}$ has size at least $t-s$ and $\mathscr{S}$ has cardinality at most $d^s$, so $\rh(\ueta_{\bB_t})\to1$ uniformly in $\ueta_{\bB_t}$ as $t\to\infty$.

Similarly, $\esmin(\cdot\giv\ueta_{\bB_t})$ is well approximated by $\esmin^0(\cdot\giv\ueta_{\bB_t})$ (within total variation distance $\oo(t)$) uniformly over all $\ueta_{\bB_t}$, since any $\si_{v_i}$ for $v_i$ not connected to $\pd\ball{t}{v_i}$ in $\ueta$ has the same distribution under the two measures independently of all the other $\si_{v_j}$, while the $\si_{v_j}$ for $v_j$ connected to $\pd\ball{t}{v_j}$ equal $1$ with probability at least $e^{Bt}/[e^{Bt}+(q-1)]$. It thus follows from \eqref{e:potts.fk} that
$\tv{
\esmin(\cdot\giv\ueta_{\bA_{s,t}})
-\esmin^0(\cdot\giv\ueta_{\bA_{s,t}})
}\le\oo(t)$ uniformly over all $\ueta_{\bA_{s,t}}$.
But under $\esmin^0$ the $\si_{v_i}$ are \emph{exactly} independent, and since $s\equiv s(t)\to\infty$ the marginal laws of the $\si_{v_i}$ resulting from the Potts--$\acr{fk}$ coupling belong to $\splxfk{\oo(t)}$, proving the result.
\epf

\pagebreak

\bibliography{potts}

\def\cprime{$'$}
\begin{thebibliography}{10}

\bibitem{Abou-Chacra}
R.~Abou-Chacra, D.~J. Thouless, and P.~W. Anderson.
\newblock A selfconsistent theory of localization.
\newblock {\em Journal of Physics C: Solid State Physics}, 6(10):1734, 1973.

\bibitem{PhysRevB.68.214403}
M.~Aizenman, R.~Sims, and S.~L. Starr.
\newblock Extended variational principle for the {S}herrington--{K}irkpatrick
  spin-glass model.
\newblock {\em Phys. Rev. B}, 68:214403, Dec 2003.

\bibitem{AizenmanWarzel}
M.~Aizenman and S.~Warzel.
\newblock The canopy graph and level statistics for random operators on trees.
\newblock {\em Mathematical Physics, Analysis and Geometry}, 9:291--333, 2006.
\newblock 10.1007/s11040-007-9018-3.

\bibitem{MR2743259}
M.~Bayati, D.~Gamarnik, and P.~Tetali.
\newblock Combinatorial approach to the interpolation method and scaling limits
  in sparse random graphs.
\newblock In {\em S{TOC}'10---{P}roceedings of the 2010 {ACM} {I}nternational
  {S}ymposium on {T}heory of {C}omputing}, pages 105--114. ACM, New York, 2010.

\bibitem{MR1873300}
I.~Benjamini and O.~Schramm.
\newblock Recurrence of distributional limits of finite planar graphs.
\newblock {\em Electron. J. Probab.}, 6(23):13 pp. (electronic), 2001.

\bibitem{Bethe}
H.~A. Bethe.
\newblock Statistical theory of superlattices.
\newblock {\em Proceedings of the Royal Society of London. Series A,
  Mathematical and Physical Sciences}, 150(871):pp. 552--575, 1935.

\bibitem{MR1757955}
M.~Biskup, C.~Borgs, J.~T. Chayes, and R.~Koteck{\'y}.
\newblock Gibbs states of graphical representations of the {P}otts model with
  external fields.
\newblock {\em J. Math. Phys.}, 41(3):1170--1210, 2000.
\newblock Probabilistic techniques in equilibrium and nonequilibrium
  statistical physics.

\bibitem{arXiv:1002.0115}
C.~Borgs, J.~Chayes, J.~Kahn, and L.~Lov\'asz.
\newblock Left and right convergence of graphs with bounded degree.
\newblock Preprint,
  \href{http://arxiv.org/abs/1002.0115v1}{\texttt{arXiv:1002.0115v1}}, 2010.

\bibitem{Chalupa}
J.~Chalupa, P.~L. Leath, and G.~R. Reich.
\newblock Bootstrap percolation on a {B}ethe lattice.
\newblock {\em Journal of Physics C: Solid State Physics}, 12(1):L31, 1979.

\bibitem{ChayesSG}
J.~T. Chayes, L.~Chayes, J.~P. Sethna, and D.~J. Thouless.
\newblock A mean field spin glass with short-range interactions.
\newblock {\em Communications in Mathematical Physics}, 106:41--89, 1986.
\newblock 10.1007/BF01210926.

\bibitem{MR534946}
V.~Chv{\'a}tal.
\newblock The tail of the hypergeometric distribution.
\newblock {\em Discrete Math.}, 25(3):285--287, 1979.

\bibitem{arXiv:1106.4714v3}
P.~Contucci, S.~Dommers, C.~Giardin\`a, and S.~Starr.
\newblock Antiferromagnetic {P}otts model on the
  {E}rd\doubleacute{o}s-{R}\'enyi random graph.
\newblock Preprint,
  \href{http://arxiv.org/abs/1106.4714v3}{\texttt{arXiv:1106.4714v3}}, 2011.

\bibitem{MR1828509}
A.~Dembo, A.~Kagan, and L.~A. Shepp.
\newblock Remarks on the maximum correlation coefficient.
\newblock {\em Bernoulli}, 7(2):343--350, 2001.

\bibitem{MR2643563}
A.~Dembo and A.~Montanari.
\newblock Gibbs measures and phase transitions on sparse random graphs.
\newblock {\em Braz. J. Probab. Stat.}, 24(2):137--211, 2010.

\bibitem{MR2650042}
A.~Dembo and A.~Montanari.
\newblock Ising models on locally tree-like graphs.
\newblock {\em Ann. Appl. Probab.}, 20(2):565--592, 2010.

\bibitem{arXiv:1110.4821}
A.~Dembo, A.~Montanari, and N.~Sun.
\newblock Factor models on locally tree-like graphs.
\newblock Preprint,
  \href{http://arxiv.org/abs/1110.4821v1}{\texttt{arXiv:1110.4821v2}}, 2011.

\bibitem{MR2733399}
S.~Dommers, C.~Giardin{\`a}, and R.~van~der Hofstad.
\newblock Ising models on power-law random graphs.
\newblock {\em J. Stat. Phys.}, 141(4):638--660, 2010.

\bibitem{GuerraRSB}
F.~Guerra.
\newblock Broken replica symmetry bounds in the mean field spin glass model.
\newblock {\em Communications in Mathematical Physics}, 233:1--12, 2003.
\newblock 10.1007/s00220-002-0773-5.

\bibitem{MR1782847}
S.~Janson, T.~{\L}uczak, and A.~Rucinski.
\newblock {\em Random graphs}.
\newblock Wiley-Interscience Series in Discrete Mathematics and Optimization.
  Wiley-Interscience, New York, 2000.

\bibitem{MR2317690}
F.~Krz{\c{a}}ka{\l}a, A.~Montanari, F.~Ricci-Tersenghi, G.~Semerjian, and
  L.~Zdeborov{\'a}.
\newblock Gibbs states and the set of solutions of random constraint
  satisfaction problems.
\newblock {\em Proc. Natl. Acad. Sci. USA}, 104(25):10318--10323 (electronic),
  2007.

\bibitem{MR1036755}
C.~McDiarmid.
\newblock On the method of bounded differences.
\newblock In {\em Surveys in combinatorics, 1989 ({N}orwich, 1989)}, volume 141
  of {\em London Math. Soc. Lecture Note Ser.}, pages 148--188. Cambridge Univ.
  Press, Cambridge, 1989.

\bibitem{MR2518205}
M.~M{\'e}zard and A.~Montanari.
\newblock {\em Information, physics, and computation}.
\newblock Oxford Graduate Texts. Oxford University Press, Oxford, 2009.

\bibitem{MezardParisiBethe}
M.~M\'ezard and G.~Parisi.
\newblock The {B}ethe lattice spin glass revisited.
\newblock {\em The European Physical Journal B - Condensed Matter and Complex
  Systems}, 20:217--233, 2001.
\newblock 10.1007/PL00011099.

\bibitem{MR1026102}
M.~M{\'e}zard, G.~Parisi, and M.~A. Virasoro.
\newblock {\em Spin glass theory and beyond}, volume~9 of {\em World Scientific
  Lecture Notes in Physics}.
\newblock World Scientific Publishing Co. Inc., Teaneck, NJ, 1987.

\bibitem{MezParZec_science}
M.~M\'ezard, G.~Parisi, and R.~Zecchina.
\newblock Analytic and algorithmic solution of random satisfiability problems.
\newblock {\em Science}, 297(5582):812--815, 2002.

\bibitem{springerlink:10.1007/s00440-010-0315-6}
A.~Montanari, E.~Mossel, and A.~Sly.
\newblock The weak limit of {I}sing models on locally tree-like graphs.
\newblock {\em Probability Theory and Related Fields}, 152:31--51, 2012.
\newblock 10.1007/s00440-010-0315-6.

\bibitem{arXiv:1203.2602}
A.~Sly and N.~Sun.
\newblock The computational hardness of counting in two-spin models on
  $d$-regular graphs.
\newblock Preprint,
  \href{http://arxiv.org/abs/1203.2602}{\texttt{arXiv:1203.2602v1}}, 2012.

\bibitem{MR2731561}
M.~Talagrand.
\newblock {\em Mean field models for spin glasses: {V}olume {I}}, volume~54 of
  {\em Ergebnisse der Mathematik und ihrer Grenzgebiete. 3. Folge. A Series of
  Modern Surveys in Mathematics}.
\newblock Springer-Verlag, Berlin, 2011.
\newblock Basic examples.

\bibitem{ThoulessSG}
D.~J. Thouless.
\newblock Spin-glass on a {B}ethe lattice.
\newblock {\em Phys. Rev. Lett.}, 56:1082--1085, Mar 1986.

\bibitem{Weiss}
P.~R. Weiss.
\newblock The application of the {B}ethe-{P}eierls method to ferromagnetism.
\newblock {\em Phys. Rev.}, 74:1493--1504, Nov 1948.

\bibitem{MR1424469}
E.~T. Whittaker and G.~N. Watson.
\newblock {\em A course of modern analysis}.
\newblock Cambridge Mathematical Library. Cambridge University Press,
  Cambridge, 1996.
\newblock Reprint of the fourth (1927) edition.

\end{thebibliography}
\bibliographystyle{abbrv}

\end{document}